\documentclass[11pt]{amsart}
\usepackage{amsfonts}
\usepackage{amssymb}
\usepackage[mathscr]{eucal}
\usepackage{graphicx}
\usepackage{tikz}
\usepackage{subfig}

\newtheorem{theorem}[]{Theorem}
\newtheorem{lemma}[theorem]{Lemma}

\theoremstyle{definition}
\newtheorem*{definition}{Definition}

\theoremstyle{remark}

\numberwithin{equation}{section}

\begin{document}

\title{The Head and Tail Conjecture for Alternating Knots}

\author{Cody Armond}
\address{Mathematics Department\\
Louisiana State University\\
Baton Rouge, Louisiana}
\email{carmond@math.lsu.edu}

\subjclass{}
\date{}

\begin{abstract}
We investigate the coefficients of the highest and lowest terms (also called the head and the tail) of the colored Jones polynomial and show that they stabilize for alternating links and for adequate links. To do this we apply techniques from skein theory.

\end{abstract}

\maketitle

\section{Introduction}

The normalized colored Jones polynomial $J_{N,L}(q)$ for a link $L$ is a sequence of Laurent polynomials in the variable $q^{1/2}$, i.e. $J_{N,L} \in \mathbb{Z}[q^{1/2},q^{-1/2}]$. This sequence is defined for $N \leq 2$ so that $J_{2,L}(q)$ is the ordinary Jones polynomial, and $J_{N,U} = 1$ where $U$ is the unknot. For links $L$ with an odd number of components (including all knots), $J_{L,N}$ is actually in $\mathbb{Z}[q,q^{-1}]$. For links with an even number of components, $q^{1/2}J_{L,N} \in \mathbb{Z}[q,q^{-1}]$.

	In \cite{dl} and \cite{dl2} Oliver Dasbach and Xiao-Song Lin showed that, up to sign, the first two coefficients and the last two coefficients of $J_{N,K}(q)$ do not depend on $N$ for alternating knots. They also showed that the third (and third to last) coefficient does not depend on $N$ so long as $N \geq 3$. This and computational data led them to believe that the $k$-th coefficient does not depend on $N$ so long as $N \geq k$. The goal of this paper is to prove this conjecture for all alternating links. It is also known that this property of $J_{N,L}(q)$ does not hold for all knots. In \cite{ad}, with Oliver Dasbach, we examined the case of the $(4,3)$ torus knot for which this property fails.

\begin{definition}
For two Laurent series $P_1(q)$ and $P_2(q)$ we define
$$P_1(q) \;\dot{=}_n\; P_2(q)$$
if after multiplying $P_1(q)$ by $\pm q^{s_1}$ and $P_2(q)$ by $\pm q^{s_2}$, $s_1$ and $s_2$ some powers, to get power series $P'_1(q)$ and $P'_2(q)$ each with positive constant term, $P'_1(q)$ and $P'_2(q)$ agree $\mod q^n$.
For example $-q^{-4} + 2 q^{-3}- 3+11 q \;\dot{=}_5\; 1-2 q+3 q^4.$
\end{definition}

Another way of phrasing the above definition is that $P_1(q) \dot{=}_n P_2(q)$ if and only if their first $n$ coefficients agree up to sign.

In \cite{ad} we defined two power series, the head and tail of the colored Jones polynomial $H_L(q)$ and $T_L(q)$.

\begin{definition}
The tail of the colored Jones polynomial of a link $L$ -- if it exists -- is a series $T_L(q)$, with
$$T_L(q) \;\dot{=}_N\; J_{L,N}(q), \text{ for all } N$$

Similarly, the head $H_L(q)$ of the colored Jones polynomial of $L$ is the tail of $J_{L,N}(q^{-1})$, which is equal to the colored Jones polynomial of the mirror image of $L$.
\end{definition}

Note that $T_L(q)$ exists if and only if $J_{L,N}(q) \;\dot{=}_N\; J_{L,N+1}(q)$ for all $N$. For example, for the first few colors $N$ the colored Jones polynomial of the knot $6_2$ multiplied by $q^{2 N^2-N-1}$ is

\medskip

{\small
\begin{tabular} {ll}
$N=2:$ & $1-2 q+2 q^2-2 q^3+2 q^4-q^5+q^6$\\
$N=3:$ & $1-2 q+4 q^3-5 q^4+6 q^6 + \dots  -q^{14}+3 q^{15}-q^{16}-q^{17}+q^{18}$\\
$N=4:$ & $1-2 q+2 q^3+q^4-4 q^5-2 q^6+\dots -2 q^{29}-3 q^{30}+3
q^{32}-q^{34}-q^{35}+q^{36}$\\
$N=5:$ & $1-2 q+2 q^3-q^4+2 q^5-6 q^6+\dots-2 q^{53}-q^{54}+4
q^{55}-q^{58}-q^{59}+q^{60}$\\
$N=6:$ & $1-2 q+2 q^3-q^4-2 q^7+q^8+\dots-3 q^{82}+3
q^{84}+q^{85}-q^{88}-q^{89}+q^{90}$\\
$N=7:$ & $1-2 q+2 q^3-q^4-2 q^6+4 q^7-3 q^8+\dots 
-q^{118}+4 q^{119}+q^{121}-q^{124}-q^{125}+q^{126}$\\
\end{tabular}
}

 This is exactly the property conjectured by Dasbach and Lin to hold for all alternating knots, and the subject of the main theorem of this paper.

\begin{theorem}
\label{thmalt}
If $L$ is an alternating link, then $J_{L,N}(q) \;\dot{=}_N\; J_{L,N+1}(q)$.
\end{theorem}

Because the mirror image of an alternating link is alternating, this theorem says that the head and the tail exists for all alternating links. Theorem \ref{thmalt} was simultaneously and independently proved by Stavros Garoufalidis and Thang Le in \cite{gl} using alternate methods.

We are also able to prove a more general theorem about $A$-adequate links.

\begin{theorem}
\label{thmad}
If $L$ is a $A$-adequate link, then $J_{L,N}(q) \;\dot{=}_N\; J_{L,N+1}(q)$.
\end{theorem}

Because all alternating links are $A$-adequate, Theorem \ref{thmad} implies Theorem \ref{thmalt}.

\subsection{Plan of paper}
In section ~\ref{Background}, we discuss definitions and basic results regarding adequate links and skein theory. In section ~\ref{main}, we present the main Lemma which is a slight generalization of a Lemma in ~\cite{ad} relating the lowest terms of the colored Jones polynomial to a certain skein theoretic graph. Finally, in section ~\ref{proof}, we present the proofs of Theorems ~\ref{thmalt} and ~\ref{thmad} using the graph discussed in section ~\ref{main}.

\subsection{Acknowledgements}
I would like to thank Oliver Dasbach for all his help and advise. I would also like to thank Pat Gilmer for teaching me all I know about skein theory.

\section{Background}
\label{Background}

\subsection{Alternating and Adequate}

\begin{figure}[htbp] %
   \centering
   \includegraphics[width=4in]{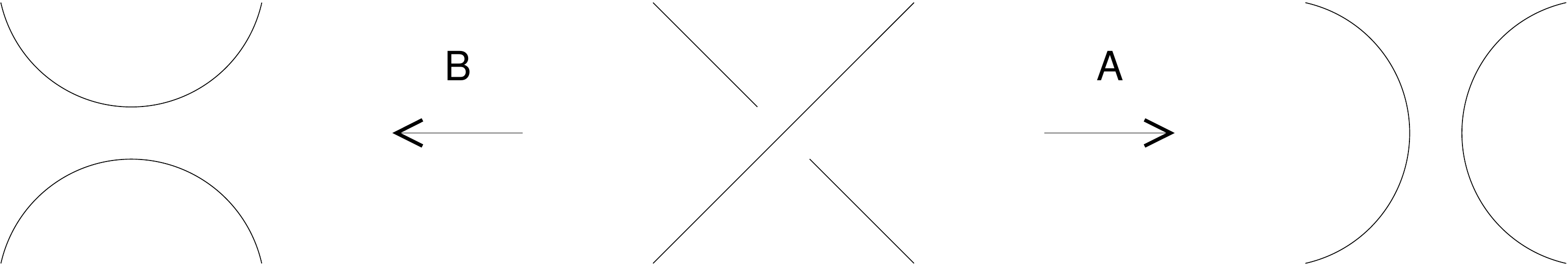} 
   \caption{A and B smoothings}
   \label{smoothings}
\end{figure}

Given a link diagram $D$ there are two ways to smooth each crossing, described in Figure \ref{smoothings}. A state of the diagram is a choice of smoothing for each crossing. Two states are particularly important when dealing with the colored Jones polynomial; they are the all-A state $S_A$ and the all-B state $S_B$. The all-A (respectively all-B) state is the state for which the A (B) smoothing is choosen for every crossing.

For a state $S$ we can build a graph $G_S$ called the state graph for $S$. The graph $G_S$ has vertices the circles in $S$, and edges the crossings in $D$. Each edge connects the two vertices corresponding to the two circles that the crossing meets.

\begin{definition}
A link diagram is A-adequate (B-adequate) if the state graph for $S_A$ ($S_B$) has no loops.

A link diagram is adequate if it is both A and B-adequate.

A link is adequate if it has an adequate diagram.
\end{definition}

The most important property of A-adequate diagrams is that the number of circles in $S_A$ is a local maximum. In other words, any state that has only a single $B$ smoothing will have one fewer circle than the all-A smoothing. Similarly for B-adequate diagrams, that the number of circles in $S_B$ is a local maximum. 

It is a well-known fact that all alternating links are adequate links. In particular, a reduced alternating diagram, that is an alternating diagram without any nugatory crossings, is an adequate diagram.

Another important fact about adequate diagrams is that parallels of A-adequate diagrams are also A-adequate diagrams. Given a diagram $D$, the $r$-th parallel of $D$ denoted $D^r$ is the diagram formed by replacing $D$ with $r$ parallel copies of $D$.

\subsection{Skein Theory} For a more detailed explanation of skein theory, see \cite{lic} or \cite{mv}

The Kauffman bracket skein module, $S(M;R,A)$, of a $3$-manifold $M$ and ring $R$ with invertible element $A$, is the free $R$-module generated by isotopy classes of framed links in $M$, modulo the submodule generated by the Kauffman relations:

\begin{center}
\begin{tabular}{c c}
$\includegraphics[width=.25in]{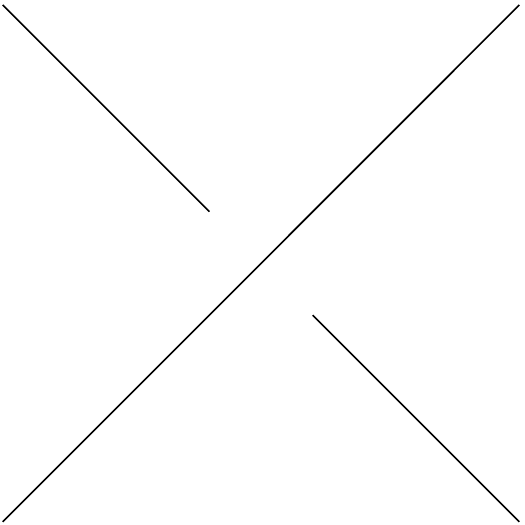} \raisebox{6pt}{$\;= A$} \includegraphics[width=.25in]{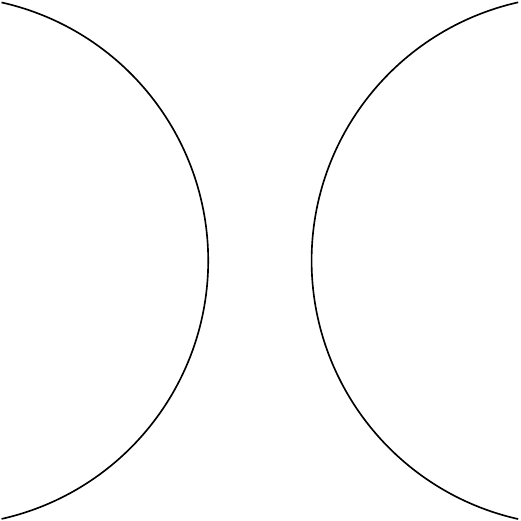} \raisebox{6pt}{$+ A^{-1}$} \includegraphics[width=.25in]{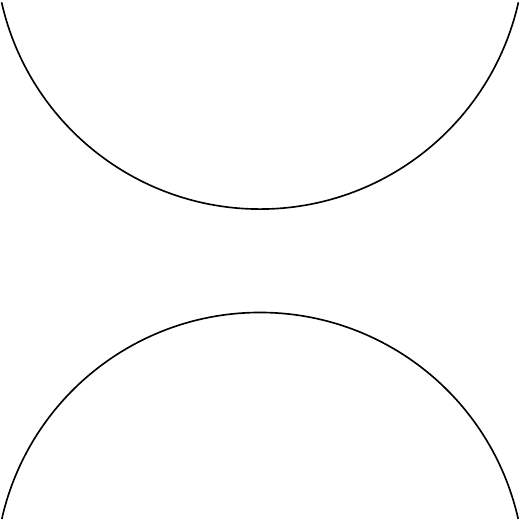}\;,\qquad\qquad$
&
$\includegraphics[width=.25in]{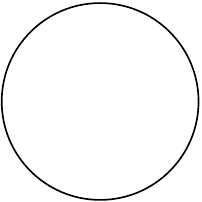} \raisebox{6pt}{$\;= -A^2 -A^{-2}$}$
\end{tabular}
\end{center}

If $M$ has designated points on the boundary, then the framed links must include arcs which meet all of the designated points.

In this paper we will take $R = \mathbb{Q}(A)$, the field of rational functions in variable $A$ with coefficients in $\mathbb{Q}$. As we are concerned with the lowest terms of a polynomial, we will need to express rational functions as Laurent series. This can always be done so that the Laurent series has a minimum degree.
\begin{definition}
Let $f\in \mathbb{Q}(A)$, define $d(f)$ to be the minimum degree of $f$ expressed as a Laurent series in $A$.
\end{definition}
Note that $d(f)$ can be calculated without referring to the Laurent series. Any rational function $f$ expressed as $\frac{P}{Q}$ where $P$ and $Q$ are both polynomials. Then $d(f) = d(P) - d(Q)$.

 We will be concerned with two particular skein modules: $S(S^3;R,A)$, which is isomorphic to $R$ under the isomorphism sending the empty link to $1$, and $S(D^3;R,A)$, where $D^3$ has $2n$ designated points on the boundary. With these designated points, $S(D^3;R,A)$ is also called the Temperley-Lieb algebra $TL_n$.

\begin{figure}[htbp]
	\centering
	\includegraphics[width=1.5in]{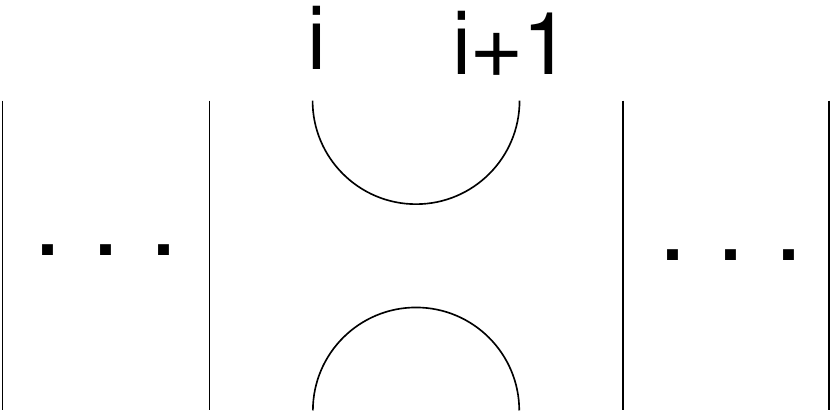}
	\caption{$h_i$}
	\label{hook}
\end{figure}

We will give an alternate explanation for the Temperley-Lieb algebra. First, consider the disk $D^2$ as a rectangle with $n$ designated points on the top and $n$ designated points on the bottom. Let $TLM_n$ be the set of all crossing-less matchings on these points, and define the product of two crossing-less matchings by placing one rectangle on top of the other and deleting any components which do not meet the boundary of the disk. With this product, $TLM_n$ is a monoid, which we shall call the Temperley-Lieb monoid. It has generators $h_i$ as in Figure \ref{hook}, and following relations:
\begin{itemize}
\item $h_ih_i=h_i$
\item $h_ih_{i\pm1}h_i = h_i$
\item $h_ih_j = h_jh_i$ if $|i-j| \geq 2$
\end{itemize}

Any element in $TL_n$ has the form $\sum_{M \in TLM_n} c_M M$, where $c_M \in \mathbb{Q}(A)$. Multiplication in $TL_n$ is slightly different from multiplication in $TLM_n$, because $h_ih_i=(-A^2-A^{-2})h_i$ in $TL_n$.

There is a special element in $TL_n$ of fundamental importance to the colored Jones polynomial, called the Jones Wentzl idempotent, denoted $f^{(n)}$. Diagramatically this element is represented by an empty box with $n$ strands coming out of it on two opposite sides. By convention an $n$ next to a strand in a diagram indicates that the strand is replaced by $n$ parallel ones.


With $$\Delta_{n}:=(-1)^{n} \frac {A^{2(n+1)}-A^{-2 (n+1)}}{A^{2}-A^{-2}}$$ and
$\Delta_{n}!:= \Delta_{n} \Delta_{n-1} \dots \Delta_{1}$ the Jones-Wenzl idempotent satisfies

\begin{center}
\begin{tabular}{c c}
$\begin{tikzpicture}[baseline=0.8cm]
\draw(0,0)--(0,1) 
node[rectangle, draw, ultra thick, fill=white]{ } --(0,2) node[right]{\scriptsize{n+1}};
\end{tikzpicture} =
\begin{tikzpicture}[baseline=0.8cm]
\draw(0,0)--(0,1) 
node[rectangle, draw, ultra thick, fill=white]{ } --(0,2) node[right]{\scriptsize{n}};
\draw(0.7,0)--(0.7,2) node[right]{\scriptsize{1}};
\end{tikzpicture} -
\left( \frac {\Delta_{n-1}} {\Delta_{n}} \right ) \, \, 
\begin{tikzpicture}[baseline=0.8cm]
\draw (0,0) node[right]{\scriptsize{n}}--(0,1) node[left]{\scriptsize{n-1}}
--(0,2) node[right]{\scriptsize{n}};
\draw[ultra thick, fill=white] (-0.2 ,0.5) rectangle (0.4, 0.65); 
\draw (0.2, 0.65) arc (180:0:0.25);
\draw (0.7,0) node[right]{1}--(0.7, 0.65);
\draw[ultra thick, fill=white] (-0.2 ,1.5) rectangle (0.4, 1.35); 
\draw (0.2, 1.35) arc (-180:0:0.25);
\draw (0.7,2) node [right]{\scriptsize{1}} --(0.7, 1.35);
\end{tikzpicture},
\, \qquad \qquad \qquad
$
&
$\begin{tikzpicture}[baseline=0.8cm]
\draw(0,0)--(0,1) 
node[rectangle, draw, ultra thick, fill=white]{ } --(0,2) node[right]{\scriptsize{1}};
\end{tikzpicture} = \, \, 
\begin{tikzpicture}[baseline=0.8cm]
\draw(0,0)--(0,2) 
node[right]{\scriptsize{1}};
\end{tikzpicture}$
\end{tabular}
\end{center}

with the properties

$$\includegraphics[height=1in]{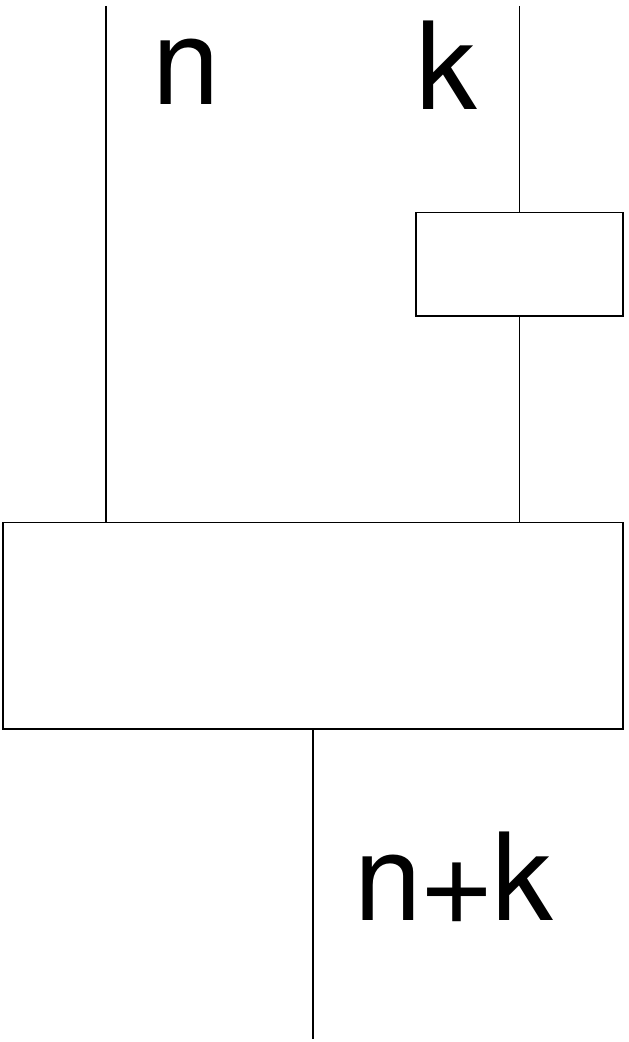} \raisebox{.5in}{$\;=\;$} \includegraphics[height=1in]{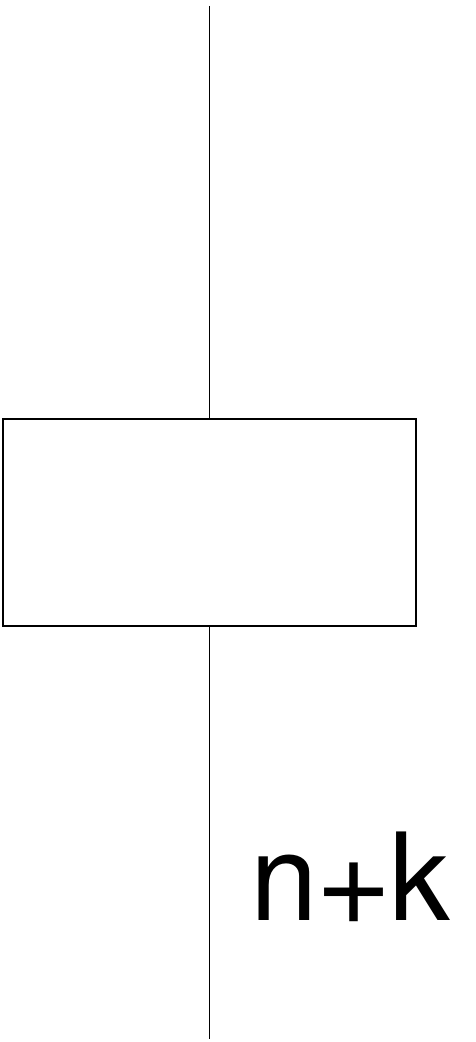}\raisebox{.5in}{,}
\, \qquad \qquad 
\includegraphics[height=1in]{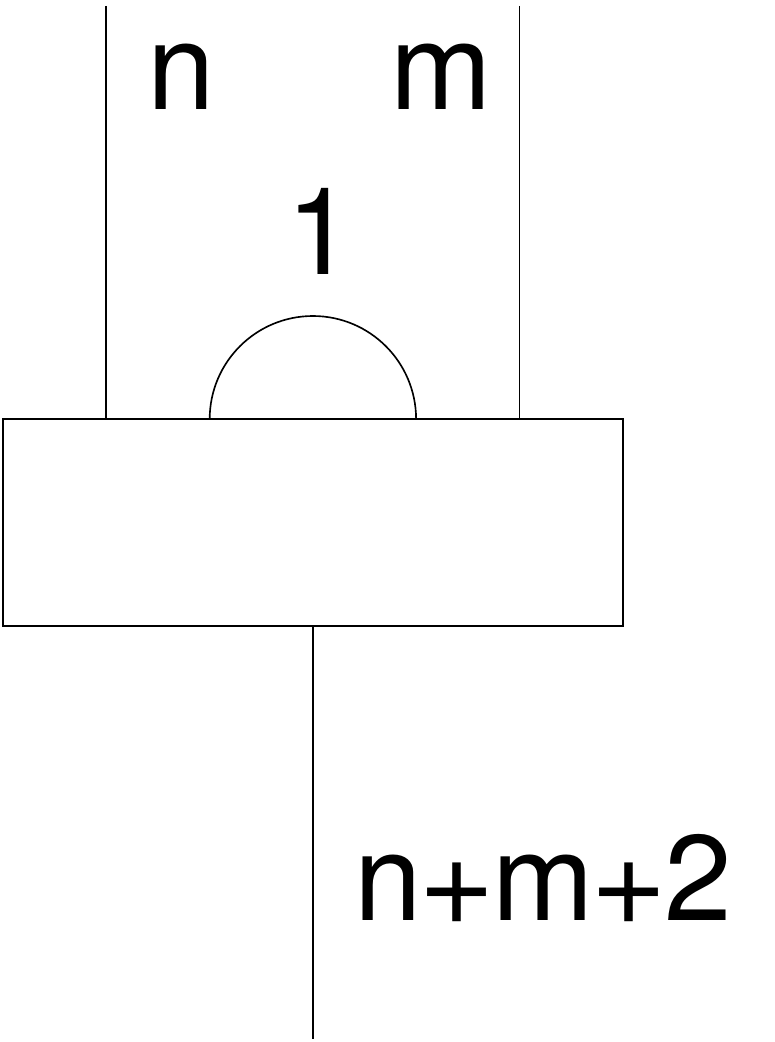} \raisebox{.5in}{$\;= 0$} 
\raisebox{.5in}{,}\, \qquad \qquad 
\includegraphics[height=1in]{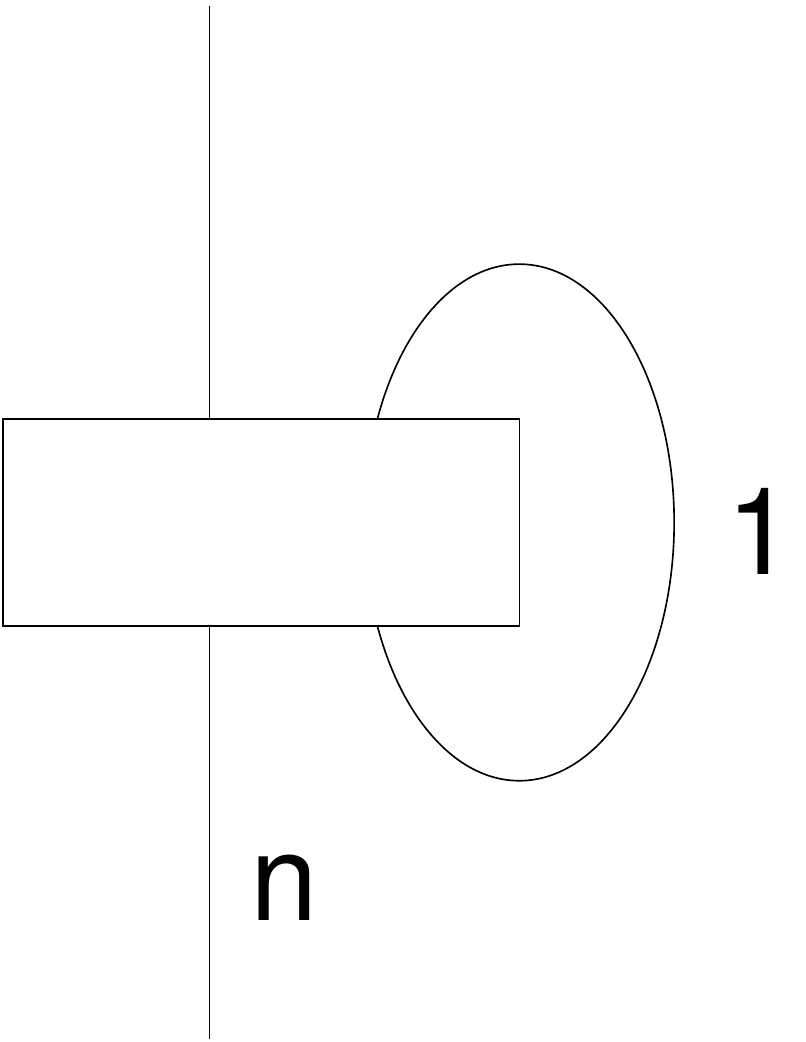} \raisebox{.45in}{$\;= \left( \frac {\Delta_{n+1}} {\Delta_{n}} \right )\;$} \includegraphics[height=1in]{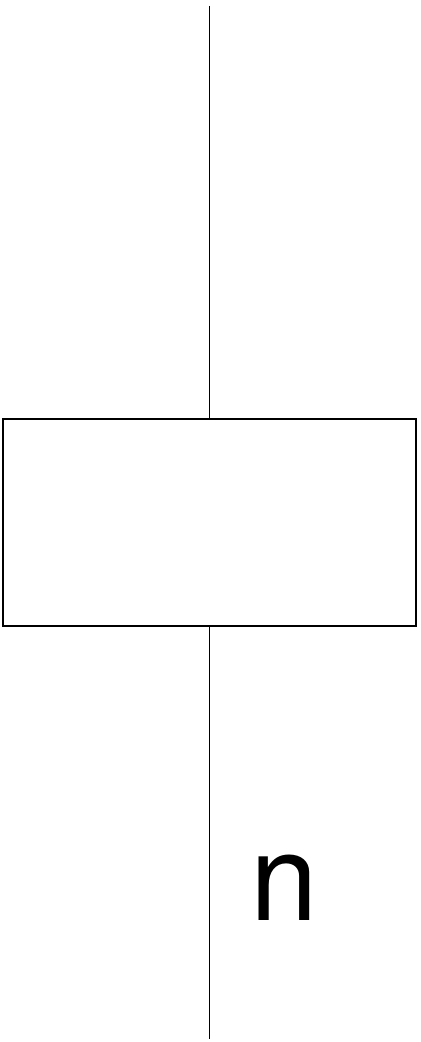}$$


If $M \in TLM_n$, define $f_M\in R$ as the coefficient of $M$ in the expansion of the Jones-Wenzl idempotent. Thus $f^{(n)} = \sum_{M\in TLM_n} f_M M$. If $e$ is the identity element of $TLM_n$, then $f_e = 1$.

\begin{lemma}
\label{idem}
If $M \in TLM_n$, then $d(f_M)$ is at least twice the minimum word length of $M$ in terms of the $h_i$'s.
\end{lemma}

\begin{proof}
This follows easily from the recursive definition of the idempotent by an inductive argument. The only issue is that terms of the form 
$\begin{tikzpicture}[baseline=0.8cm]
\draw (0,0) node[right]{\small{n}}--(0,1) node[left]{\scriptsize{n-1}}
--(0,2) node[right]{\small{n}};
\draw[ultra thick, fill=white] (-0.2 ,0.5) rectangle (0.4, 0.65); 
\draw (0.2, 0.65) arc (180:0:0.25);
\draw (0.7,0) node[right]{1}--(0.7, 0.65);
\draw[ultra thick, fill=white] (-0.2 ,1.5) rectangle (0.4, 1.35); 
\draw (0.2, 1.35) arc (-180:0:0.25);
\draw (0.7,2) node [right]{\small{1}} --(0.7, 1.35);
\end{tikzpicture}$ may have a circle which needs to be removed. In this situation, the minimum degree of the coefficient is reduced by two, but the number of generators used is also reduced by one.

\end{proof}

Using Lemma \ref{idem} we can find a lower bound for the minimum degree of any element of $S(S^3;R,A)$ which contains the Jones-Wenzl idempotent. Before we do this, consider a crossing-less diagram $S$ in the plane consisting of arcs connecting Jones-Wenzl idempotents. We will define what it means for such a diagram to be adequate in much the same way that a knot diagram can be A or B- adequate.

\begin{figure}[htbp]
	\centering
	\subfloat[An adequate diagram]{
	\includegraphics[width=2in]{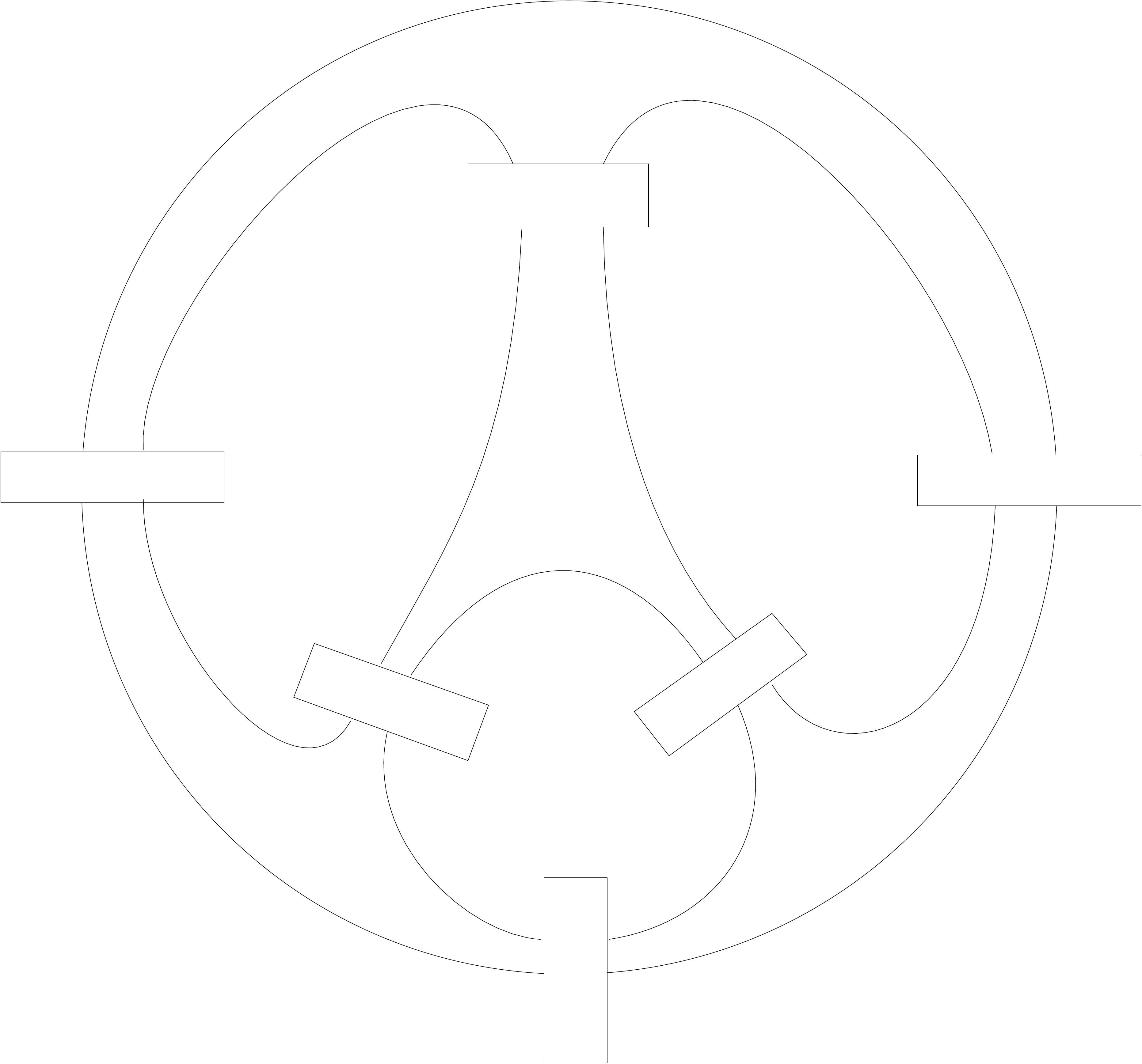}
	\label{fig:adequate}}
	\qquad
	\subfloat[An inadequate diagram]{
	\includegraphics[width=2in]{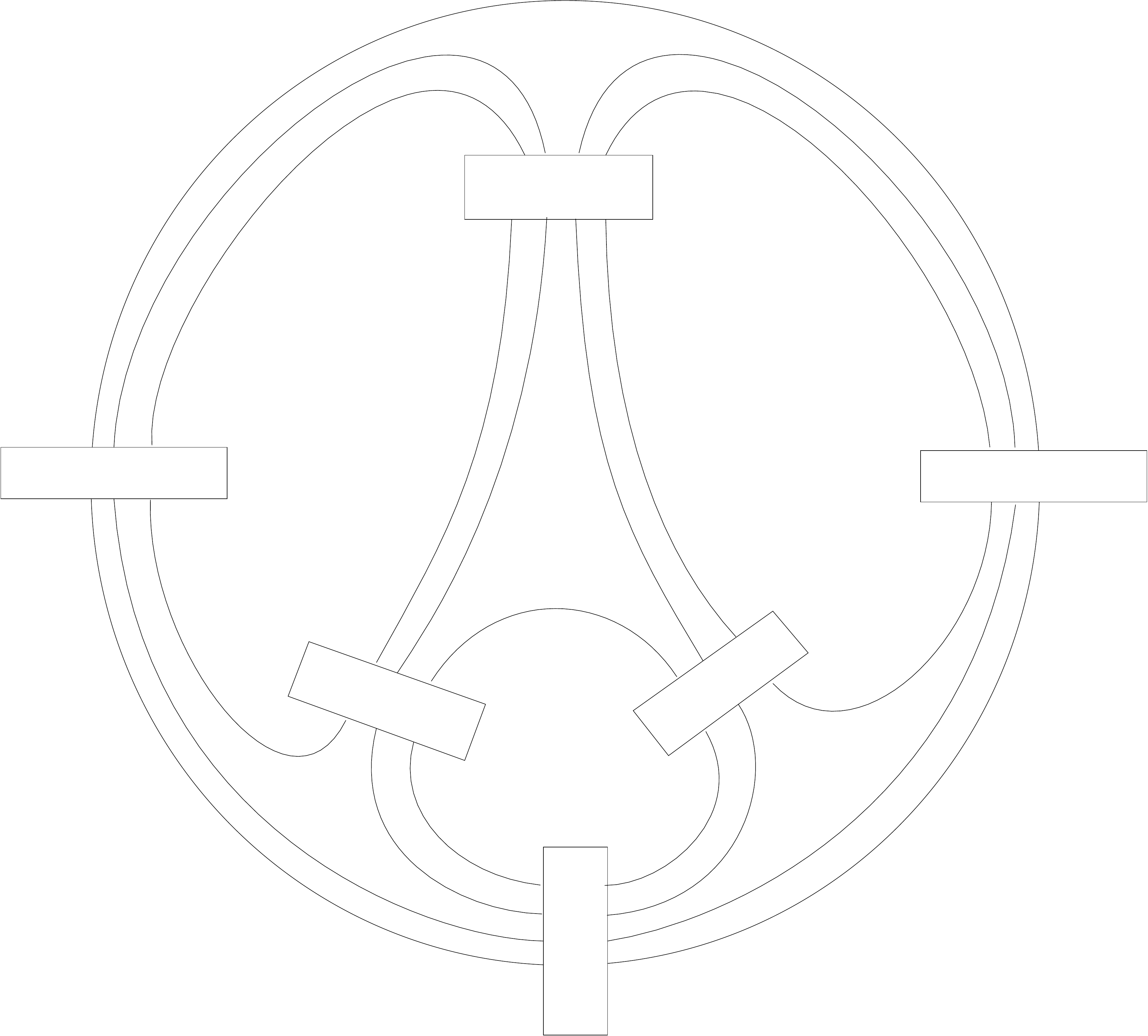}
	\label{fig:inadequate}}
	\caption{Example of adequate and inadequate diagrams}
	\label{fig:ExamplePiplup}
\end{figure}

Construct a crossing-less diagram $\bar{S}$ from $S$ by replacing each of the Jones-Wenzl idempotents in $S$ by the identity of $TL_n$. Thus $\bar{S}$ is a collection of circles with no crossings. Consider the regions in $\bar{S}$ where the idempotents had previously been. $S$ is adequate if no circle in $\bar{S}$ passes through any one of these regions more than once. Figure \ref{fig:adequate} shows an example of a diagram that is adequate and Figure \ref{fig:inadequate} shows an example of a diagram that is not adequate. In both Figures every arc is labelled $1$.

If $S$ is adequate, then the number of circles in $\bar{S}$ is a local maximum, in the sense that if the idempotents in $S$ are replaced by other elements of $TLM_n$ such that there is exactly one hook total in all of the replacements, then the number of circles in this diagram is one less than the number of circles in $\bar{S}$.

If the diagram $S$ happens to have crossings in it, we can still construct the diagram $\bar{S}$, which is now a link diagram. Denote $D(S) := d(\bar{S})$.

\begin{lemma}
\label{replace}
If $S \in S(S^3;R,A)$ is expressed as a single diagram containing the Jones-Wenzl idempotent, then $d(S) \geq D(S)$.

If the diagram for $S$ is a crossing-less adequate diagram, then $d(S) = D(S)$.
\end{lemma}

\begin{proof}

First suppose that the diagram $S$ has no crossings. We can get an expansion of $S$ by expanding each of the idempotents that appear in the diagram. Consider a single term $T_1$ in this expansion. Unless all of the idempotents have been replaced by the identity in this term, then there will be a hook somewhere in the diagram. By removing a single hook, we get a different term $T_2$ in the expansion. The number of circles in $T_1$ differs from the number of circles in $T_2$ by exactly one. Also there are fewer hooks in $T_2$, so by Lemma \ref{idem} and the fact that removing a circle results in multiplying by $-A^2-A^{-2}$, the minimum degree of $T_1$ is at least as large as the minimum degree of $T_2$. This tells us that the lowest degree amongst terms in the expansion of $S$ is the degree of the term with the idempotents replaced by the identity, $\bar{S}$.

If $S$ is adequate, then for any term $T_1$ with only a single hook, $T_2$ will be $\bar{S}$, and thus $T_2$ will have one more circle than $T_1$. Therefore, $d(T_1)>d(\bar{S})$. This tells us that any term $T$ in this expansion will have $d(T_1)>d(\bar{S})$, and thus $d(S)=d(\bar{S})=D(S)$.

If there were crossings in $S$, then we can get an expansion of $S$ by expanding the idempotents that appear in $S$ and summing over all possible smoothings of the crossings. If we expand over the smoothings first, we get a collection of terms each of which is a coefficient times a crossingless diagram with idempotents. We can apply the previous argument to say that the minimum degree of each term is the minimum degree of that term with the idempotents replaced with the identity. Now consider $\bar{S}$. By expanding $\bar{S}$ by summing over all possible smoothings of $\bar{S}$, we get the same sum as before. Thus the minimum degree of $\bar{S}$ agrees with the minimum degree of $S$.
\end{proof}

We can use trivalent graphs to express elements in a skein module using the following correspondence:

$$\includegraphics[height=1in]{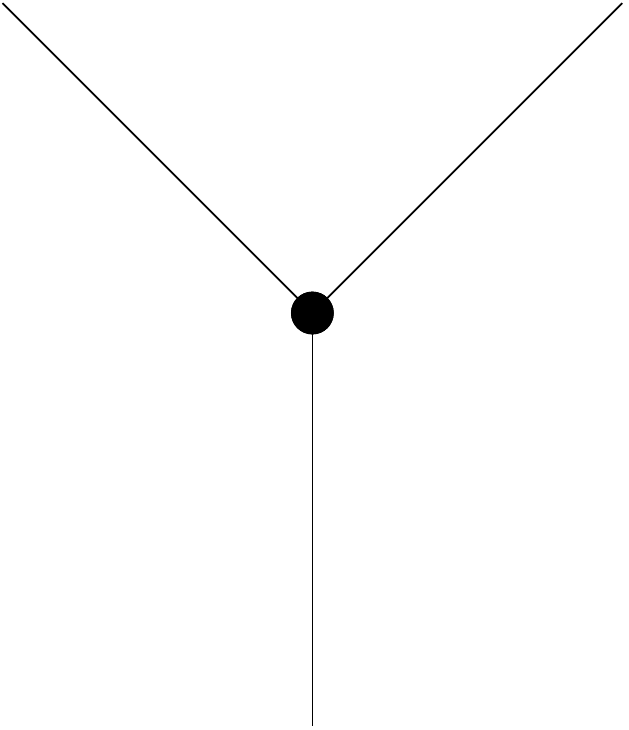} \raisebox{.5in}{$\;=\;$} \includegraphics[height=1in]{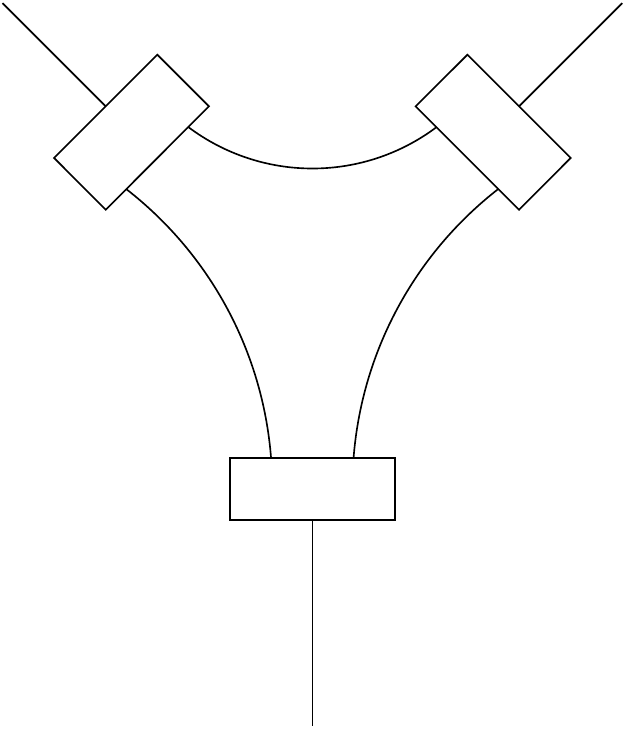}$$

Fusion is given by

$$\begin{tikzpicture}[baseline=0.8cm]
 \draw (0.7,0) .. controls (0.2 ,0.6) and (0.2,1.4) .. (0.7,2) node[right]{\scriptsize{b}};
 \draw (0.35,1) node[shape=rectangle, draw, ultra thick, fill=white]{ };
\draw (-0.7,0).. controls (-.2,0.6) and (-0.2,1.4)  .. (-0.7,2) node[left]{\scriptsize{a}};
 \draw (-0.35,1) node[rectangle, draw, ultra thick, fill=white]{ };
\end{tikzpicture}
=  \sum_{c} \frac {\Delta_c} {\theta(a,b,c)} 
\begin{tikzpicture}[baseline=0.8cm]
 \draw (1,0) node[right]{\scriptsize{b}} --(0.5,0.5);
 \draw (0,0) node[left]{\scriptsize{a}} --(0.5,0.5);
 \draw (0.5,0.5)--(0.5,1) node[right]{c}-- (0.5,1.5);
\draw (0.5,1.5)--(0,2) node[left] {\scriptsize{a}};
\draw (0.5,1.5)--(1,2) node[right] {\scriptsize{b}};
   \fill (0.5,0.5) circle (0.07);
   \fill (0.5,1.5) circle (0.07);
\end{tikzpicture} 
$$

where the sum is over all $c$ such that:
\begin{enumerate}
\item $a+b+c$ is even
\item $|a -b | \leq c \leq a+b.$
\end{enumerate}
 
\bigskip To define $\theta(a,b,c)$ let $a,b$ and $c$ related as above and  $x, y$ and $z$ be defined by $a=y+z, b=z+x$ and $c=x+y$ then

$$\theta(a,b,c):= \, \begin{tikzpicture}[baseline=0.8cm, rounded corners=2mm]
\draw (0,1) ellipse (0.6 and 1) ;
\draw (0.7,0.96) node[right]{\scriptsize{c}} ;
\draw (-0.7,0.96) node[left]{\scriptsize{a}};
\draw (0,0)--(0,1) node[right]{\scriptsize{b}}--(0,2);
   \fill (0,0) circle (0.07);
   \fill (0,2) circle (0.07);
\end{tikzpicture}$$

and one can show that

$$\theta(a,b,c)= \frac { \Delta_{x+y+z}! \Delta_{x-1}! \Delta_{y-1}! \Delta_{z-1}! } {\Delta_{y+z-1}! \Delta_{z+x-1}! \Delta_{x+y-1}!}.$$

Furthermore one has:

$$\begin{tikzpicture}[baseline=.5cm, rounded corners=2mm]
    \draw(0,0)-- (-.5,.5) -- (-.1,0.9);
    \draw (.1,1.1)--(.5,1.5) node[right]{\scriptsize{b}};
      \draw (0,0) -- (.5,.5)--(-.5,1.5) node[left]{\scriptsize{a}};
     \draw (0,0)--(0,-0.8) node[right]{\scriptsize{c}};
     \fill (0,0) circle (0.07);
\end{tikzpicture}  = (-1)^{\frac{a+b-c} 2} A^{a+b-c+\frac{a^2+b^2-c^2}{2}}
\begin{tikzpicture}[baseline=.5cm, rounded corners=2mm]
       \draw (0,0)--(0,-0.8) node[right]{\scriptsize{c}}; 
     \draw (0,0)--(0.5,1.5) node[right]{\scriptsize{b}};
     \draw (0,0)--(-0.5,1.5) node[left]{\scriptsize{a}};
      \fill (0,0) circle (0.07);
      \end{tikzpicture}$$

We are only interested in the list of coefficients of the colored Jones polynomial. In particular we consider polynomials up to powers of their variable. Up to a factor of $\pm A^s$ for some power $s$ that depends on the writhe of the link diagram the (unreduced) colored Jones polynomial $\tilde J_{n,L}(A)$ of a link $L$ can be defined as the value of the skein relation applied to the link were every component is decorated by an $n$ together with the Jones-Wenzl idempotent.  Recall that $A^{-4}=q$. To obtain the reduced colored Jones polynomial we have to divide $\tilde J_{n,L}(A)$ by its value on the unknot. Thus 
$$J_{n+1,L}(q):=\left.\frac{\tilde J_{n,L}(A)} {\Delta_{n}}\right|_{A=q^{-1/4}}.$$

\section{The Main Lemma}
\label{main}


In this section we will relate the tail of the colored Jones polynomial to a certain trivalent graph viewed as an element of the Kauffman bracket skein module of $\mathbb{R}^3$. This construction was used in \cite{ad} to prove interesting properties of the head and the tail of the colored Jones polynomial.

Given an B-adequate diagram $D$ of a link $L$, consider a negative twist region. Apply the identities of Section \ref{Background} to get the equation:

$$\begin{tikzpicture}[scale=0.6, baseline=0.8cm, rounded corners=2mm]
    \draw (-.5,2) -- (-0.1, 2.4);
    \draw (0.1,2.6) -- (0.5,3);
    \draw (0.5, 2) -- (-0.5,3);
    \draw[dashed] (0, 1.5) -- (0,1.75) node[left]{\small{m}} --(0, 2);
    \draw (-.5, -0.5) node[left]{\tiny{n}} -- (-0.1,-0.1);
    \draw (.5,-0.5) node[right]{\tiny{n}}-- (0,0)-- (-.5,.5) -- (-.1,0.9);
    \draw (.1,1.1)--(.5,1.5);
      \draw (0.1,0.1) -- (.5,.5)--(-.5,1.5) ;
\end{tikzpicture} \, = \sum_{j=0}^n (\gamma(n,n,2j))^m \frac{\Delta_{2j}}{\theta(n,n,2j)} 
\begin{tikzpicture}[scale=0.6, baseline=0.8cm]
    \draw (-.5,3) node[left]{\tiny{n}}-- (0, 2.5);
    \draw (0.5,3) node[right]{\tiny{n}}-- (0,2.5);
    \draw (0, 2.5) --(0,1.25) node[right]{\tiny{2j}}-- (0,0);
    \draw (0,0)--(-.5,-.5) node[left]{\tiny{n}};
    \draw (0,0)--(.5,-.5) node[right]{\tiny{n}};
\end{tikzpicture}
.$$

Here $\gamma(a,b,c):=(-1)^{\frac {a+b-c} 2} A^{a+b-c+ \frac{a^2+b^2-c^2} 2}.$

If we apply this equation to every maximal negative twist region, then we get an embedded trivalent graph called $\Gamma$. We get a colored graph $\Gamma_{n,(j_1,\ldots,j_k)}$ where $k$ is the number of maximal negative twist regions and $0\leq j_i \leq n$ by coloring the edges coming from the $i$-th twist region by $2j_i$ and coloring all of the other edges by $n$. From the previous equation, it is clear that we get
$$\tilde J_{n,L} \dot{=} \sum_{j_1,\ldots,j_k=0}^n \prod_{i=1}^k (\gamma(n,n,2j_i))^m \prod_{i=1}^k \frac{\Delta_{2j_i}}{\theta(n,n,2j_i)} \Gamma_{n,(j_1,\ldots,j_k)}$$

The following Theorem is a useful tool to find properties of the head and tail of the colored Jones polynomial. In this paper, it will be used to prove the existence of the head and tail for all adequate links:

\begin{theorem}\label{skein}
If $D$ is a B-adequate diagram of the link $L$, and $\Gamma_{n,(n,\ldots,n)}$ is the corresponding graph, then
$$\tilde J_{n,L} \dot{=}_{4(n+1)} \Gamma_{n,(n,\ldots,n)}$$
\end{theorem}

This Theorem was proved for the case when $D$ is a reduced alternating diagram in \cite{ad} as Theorem 4.3. The proof given there extends easily to B-adequate diagrams. We will present the proof again here with the modifications. In later sections, we shall denote $\Gamma_{n} := \Gamma_{n,(n,\ldots,n)}$.

For a rational function $R$, let $d(R)$ be the minimum degree of $R$ considered as a power series. The theorem will now follow from the following three lemmas.

\begin{lemma}
\label{gamma}
\begin{eqnarray*}
d(\gamma(n,n,2n)) &=& d(\gamma(n,n,2(n-1))) - 4n\\
d(\gamma(n,n,2j)) &\leq& d(\gamma(n,n,2(j-1)))
\end{eqnarray*}
\end{lemma}

\begin{lemma}
\label{fuse}
$$d\left(\frac{\Delta_{2j}}{\theta(n,n,2j)}\right) = d\left(\frac{\Delta_{2(j-1)}}{\theta(n,n,2(j-1))}\right) - 2$$
\end{lemma}

\begin{lemma}
\label{graph}
If $\Gamma$ is the graph coming from a B-adequate diagram, then
\begin{eqnarray*}
D(\Gamma_{n,(j_1,\ldots,j_{(i-1)},j_i,j_{i+1},\ldots,j_k)}) &=& D(\Gamma_{n,(j_1,\ldots,j_{(i-1)},j_i-1,j_{i+1},\ldots,j_k)}) \pm 2\\
d(\Gamma_{n,(n,\ldots,n,\ldots,n)}) &=& D(\Gamma_{n,(n,\ldots,n-1,\ldots,n)}) - 2
\end{eqnarray*}
\end{lemma}

\begin{proof}[Proof of Lemma \ref{gamma}]
\begin{eqnarray*}
\gamma(n,n,2j) &=& \pm A^{n+n-2j+\frac{n^2+n^2-(2j)^2}{2}}\\
&=& \pm A^{2n-2j+n^2-2j^2}
\end{eqnarray*}
Clearly $d(\gamma(n,n,2j))$ increases as $j$ decreases. Furthermore:
\begin{eqnarray*}
d(\gamma(n,n,2n)) &=& -n^2\\
d(\gamma(n,n,2(n-1))) &=& 2n-2(n-1)+n^2-2(n-1)^2\\
&=&-n^2 + 4n\\
\end{eqnarray*}

\end{proof}

\begin{proof}[Proof of Lemma \ref{fuse}]

To calculate $\theta(n,n,2j)$ note that in the previous formula for $\theta$ we get $x= j$, $y=j$, and $z=n-j$. Using this and the fact that $d(\Delta_n) = -2n$, we get:

\begin{eqnarray*}
d\left(\frac{\Delta_{2j}}{\theta(n,n,2j)}\right) &=& d\left(\frac{\Delta_{2j}\Delta_{n-1}!\Delta_{n-1}!\Delta_{2j-1}!}{\Delta_{n+j}!\Delta_{j-1}!\Delta_{j-1}!\Delta_{n-j-1}!}\right)\\
&=& d\left(\frac{\Delta_{2j}\Delta_{2j-1}\Delta_{n-j}}{\Delta_{n+j}\Delta_{j-1}\Delta_{j-1}}\right) + d\left(\frac{\Delta_{2(j-1)}!\Delta_{n-1}!\Delta_{n-1}!}{\Delta_{n+j-1}!\Delta_{j-2}!\Delta_{j-2}!\Delta_{n-j}!}\right)\\
&=& -4j -2(2j-1) - 2(n-j)+4(j-1) +2(n+j) + d\left(\frac{\Delta_{2(j-1)}}{\theta(n,n,2(j-1))}\right)\\
&=& -2 + d\left(\frac{\Delta_{2(j-1)}}{\theta(n,n,2(j-1))}\right)
\end{eqnarray*}

\end{proof}





\begin{proof}[Proof of Lemma \ref{graph}]

Consider the graph $\Gamma_{n,(j_1,\ldots,j_k)}$ viewed as an element in $S(S^3;\mathbb{Q}(A),A)$. We must compare $D(\Gamma_{n,(j_1,\ldots,j_{(i-1)},j_i,j_{i+1},\ldots,j_k)})$ with $D(\Gamma_{n,(j_1,\ldots,j_{(i-1)},j_i-1,j_{i+1},\ldots,j_k)})$. Recall that $D(S)$ is $-2$ times the number of circles in $\bar{S}$, where $\bar{S}$ is obtained from $S$ by replacing the idempotents in the diagram by the identity in $TL_m$. For $\Gamma_{n,(j_1,\ldots,j_{(i-1)},j_i,j_{i+1},\ldots,j_k)}$ and $\Gamma_{n,(j_1,\ldots,j_{(i-1)},j_i-1,j_{i+1},\ldots,j_k)}$, the number of circles in each diagram differ by $1$. Thus $D(\Gamma_{n,(j_1,\ldots,j_{(i-1)},j_i,j_{i+1},\ldots,j_k)}) = D(\Gamma_{n,(j_1,\ldots,j_{(i-1)},j_i-1,j_{i+1},\ldots,j_k)}) \pm 2$.

For $\Gamma_{n,(n,\ldots,n)}$, replacing the idempotents with the identitiy results in the all B-smoothing of the diagram $D^n$. Since $D$ is a B-adequate diagram, so is $D^n$. For $\Gamma_{n,(n,\ldots,n-1,\ldots,n)}$, the replacement results in a smoothing of $D^n$ with exactly one A smoothing. Thus the result of the replacement for $\Gamma_{n,(n,\ldots,n)}$ will have one more circle than the result of the replacement for $\Gamma_{n,(n,\ldots,n-1,\ldots,n)}$, which give us $D(\Gamma_{n,(n,\ldots,n,\ldots,n)}) = D(\Gamma_{n,(n,\ldots,n-1,\ldots,n)}) - 2$.

Finally, since $D^n$ is B-adequate, $\Gamma_{n,(n,\ldots,n)}$ is adequate. Thus by Lemma \ref{replace}, $$d(\Gamma_{n,(n,\ldots,n,\ldots,n)}) = D(\Gamma_{n,(n,\ldots,n,\ldots,n)}).$$

\end{proof}

\section{Proof of the Main Theorem}
\label{proof}
Using Theorem \ref{skein}, the main Theorem can be rewritten as follows:

\begin{theorem}
\label{restated}
If $D$ is a B-adequate diagram for a link $L$ and $\Gamma_n$ its corresponding graph, then

$$\Gamma_n \dot{=}_{4(n+1)} \Gamma_{n+1}$$

\end{theorem}

\begin{proof}

We will first prove Theorem \ref{restated} in the case of $D$ being a reduced alternating diagram, and then we will show how the proof can be modified to apply to any B-adequate diagram in general.

Interpretting $\Gamma_n$ as a skein element, we may use the following simplification:

$$
\begin{tikzpicture}[baseline=2.8cm]
\draw (1,1.6) node[right]{\scriptsize{2n}}-- (1,3);
\draw (1,3)--(0,4) node[left] {\scriptsize{n}};
\draw (1,3)--(2,4) node[right] {\scriptsize{n}};
    \fill (1,3) circle (0.07);
\end{tikzpicture}
=
\begin{tikzpicture}[baseline=2.8cm]
\draw (1.7,3.5)--(1.7,4) node[right] {\scriptsize{n}};
\draw (0.3,3.5)--(0.3,4) node[left] {\scriptsize{n}};
\draw[ultra thick, fill=white] (1.4,3.35) rectangle (2.0, 3.5);
\draw[ultra thick, fill=white] (0,3.35) rectangle (0.6, 3.5);
\draw (0.45,3.35) .. controls (0.8,2.9) and (1.2,2.9) .. (1.55, 3.35);
\draw (1, 3.2) node {\scriptsize{0}};
\draw (1,1.6) node[right]{\scriptsize{2n}}-- (1,2.1);
\draw[ultra thick, fill=white] (0.7,2.1) rectangle (1.3, 2.25);
\draw (0.85,2.25)-- (0.15, 3.35)  ;
\draw (0.45, 2.8) node[left] {\scriptsize{n}};
\draw(1.15,2.25)--(1.85, 3.35);
\draw (1.5,2.8) node[right] {\scriptsize{n}};
\end{tikzpicture} 
=
\begin{tikzpicture}[baseline=2.8cm]
\draw (1.7,3.5)--(1.7,4) node[right] {\scriptsize{n}};
\draw (0.3,3.5)--(0.3,4) node[left] {\scriptsize{n}};
\draw[ultra thick, fill=white] (1.4,3.35) rectangle (2.0, 3.5);
\draw[ultra thick, fill=white] (0,3.35) rectangle (0.6, 3.5);
\draw (1,1.6) node[right]{\scriptsize{2n}}-- (1,2.1);
\draw[ultra thick, fill=white] (0.7,2.1) rectangle (1.3, 2.25);
\draw (0.85,2.25)-- (0.15, 3.35)  ;
\draw (0.45, 2.8) node[left] {\scriptsize{n}};
\draw(1.15,2.25)--(1.85, 3.35);
\draw (1.5,2.8) node[right] {\scriptsize{n}};
\end{tikzpicture} 
=
\begin{tikzpicture}[baseline=2.8cm]
\draw (1,1.6) node[right]{\scriptsize{2n}}-- (1,3);
\draw (1,3)--(0,4) node[left] {\scriptsize{n}};
\draw (1,3)--(2,4) node[right] {\scriptsize{n}};
    \draw[ultra thick, fill=white] (0.7,2.9) rectangle (1.3, 3.1);
\end{tikzpicture}
$$

In general, $\Gamma_n$ will reduce to a collection of circles coming from the all-B smoothing $S_B$, ``fused'' together with the Jones-Wenzl Idempotent colored $2n$ for each maximal negative twist region. We will call this reduced form $S^{(n)}_B$.

\begin{figure}[htbp] %
   \centering
   \subfloat[The knot $6_2$]{
   \includegraphics[width=1.5in]{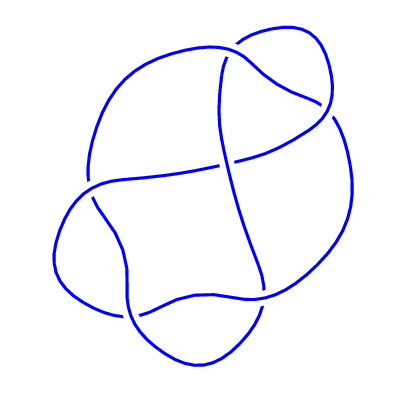}}
   \raisebox{42pt}{$\longrightarrow$}
   \subfloat[$\Gamma_n$]{
   \includegraphics[width=1in]{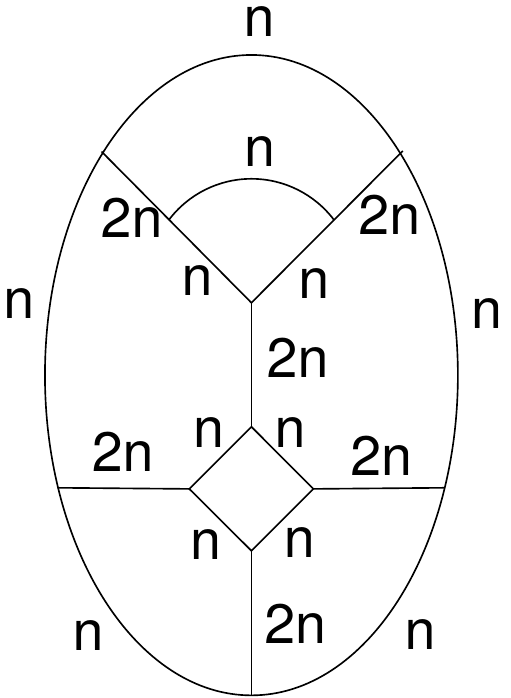}} 
   \raisebox{42pt}{$\longrightarrow$}
   \subfloat[$S^{(n)}_B$]{
   
   \includegraphics[width=1in]{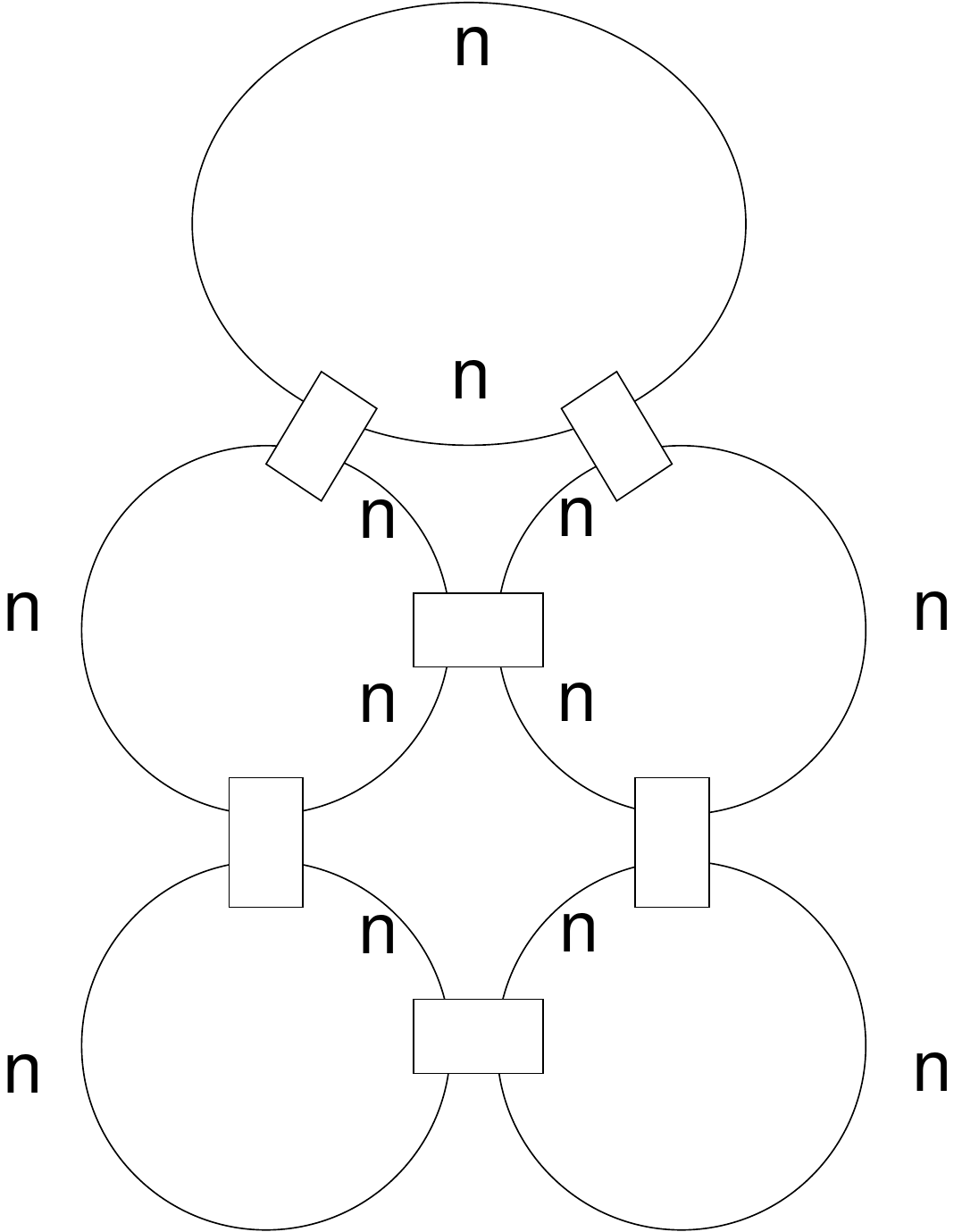}}
   \caption{Example of the knot $6_2$ along with $\Gamma_n$ and $S^{(n)}_B$}
   \label{exsB}
\end{figure}

We would now like to consider $S^{(n+1)}_B$ and show that we can reduce it to $S^{(n)}_B$ without affecting the lowest $4(n+1)$ terms. To do this we will first show a local relation which we will be able to use repeatedly.

\begin{lemma}
\label{skeinreduc}
$$\raisebox{3pt}{\includegraphics[width=1in]{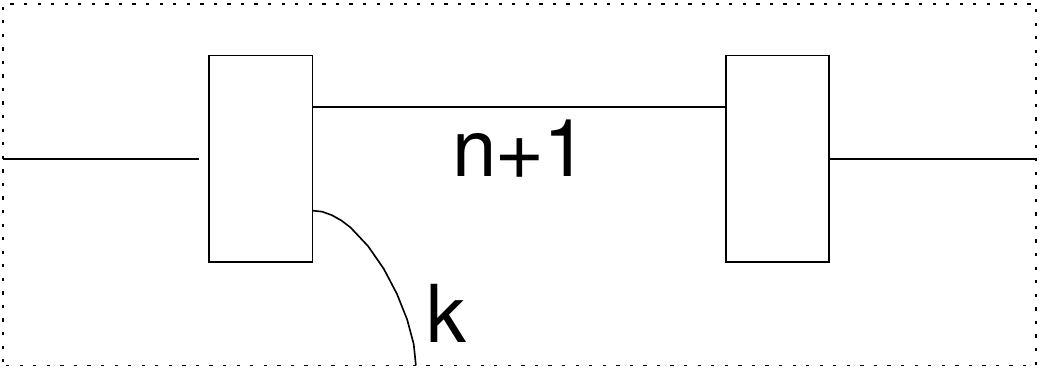}} \raisebox{12pt}{$\; = \;$} \includegraphics[width=1in]{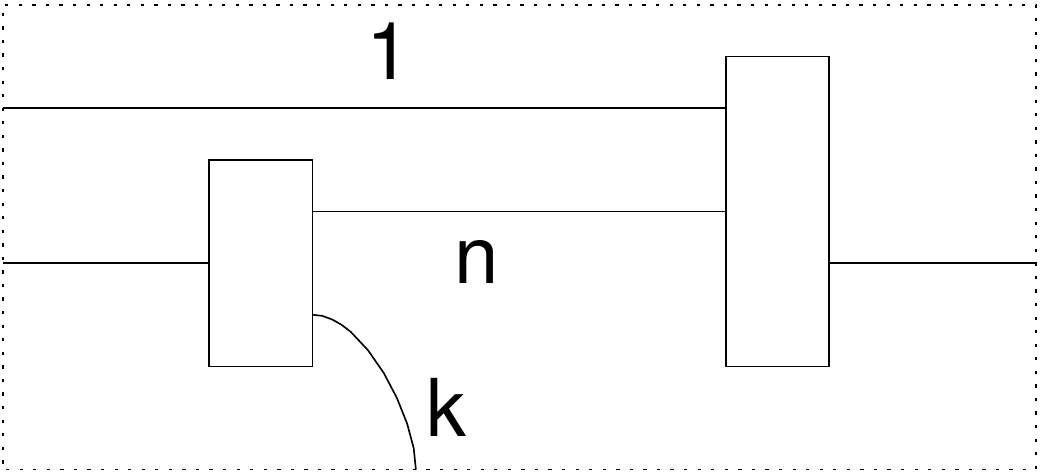} \raisebox{12pt}{$\; - \left(\frac{\Delta_{k-1}}{\Delta_{n+k}}\right)\;$} \raisebox{3pt}{\includegraphics[width=1in]{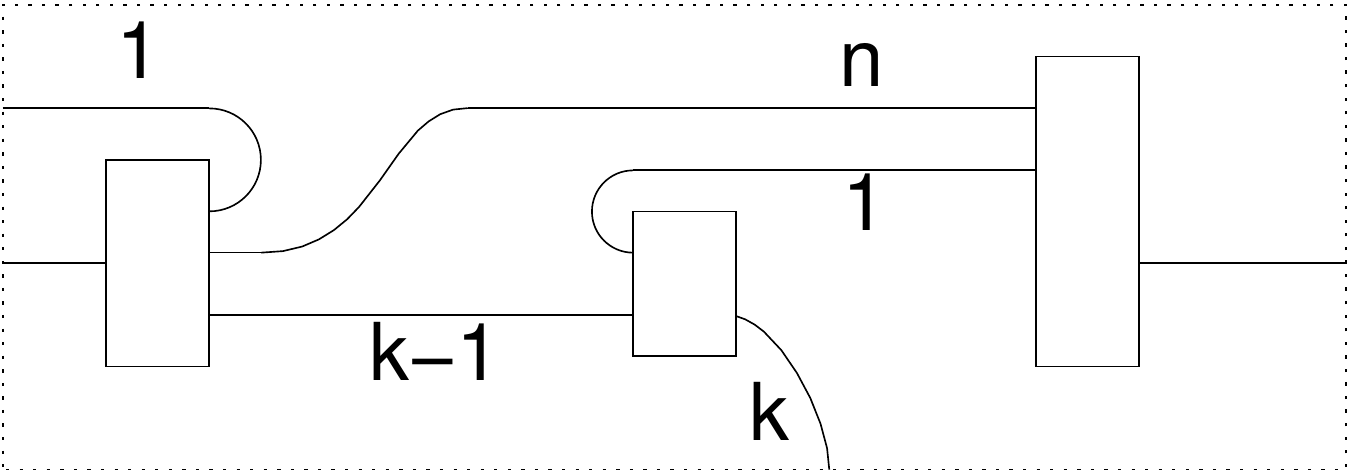}}$$
\end{lemma}

\begin{proof}
Using the recursive formula for the Jones-Wenzl idempotent on the left of the left hand side of the identity we get:
$$\raisebox{3pt}{\includegraphics[width=1in]{Skeinlemma1v2.pdf}} \raisebox{12pt}{$\; = \;$} \includegraphics[width=1in]{Skeinlemma2v2.pdf} \raisebox{12pt}{$\; - \left(\frac{\Delta_{n+k-1}}{\Delta_{n+k}}\right)\;$} \raisebox{3pt}{\includegraphics[width=1in]{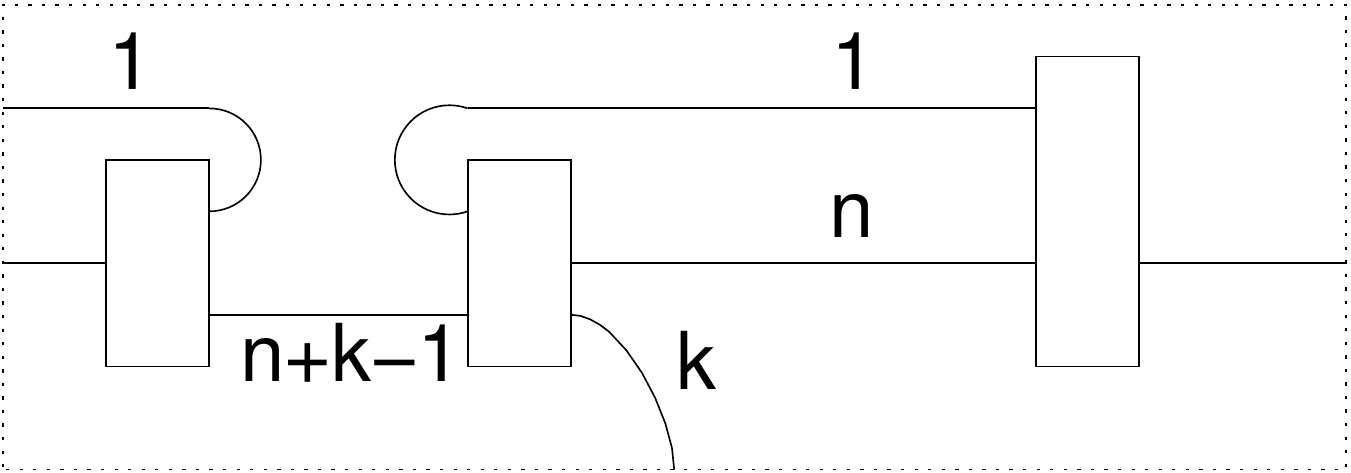}}$$
Applying the recursive formula again on the middle idempotent of the right most diagram we get:
\begin{eqnarray*}
\includegraphics[width=1in]{Skeinlemma3v2.pdf} &\raisebox{8pt}{$=$}& \includegraphics[width=1in]{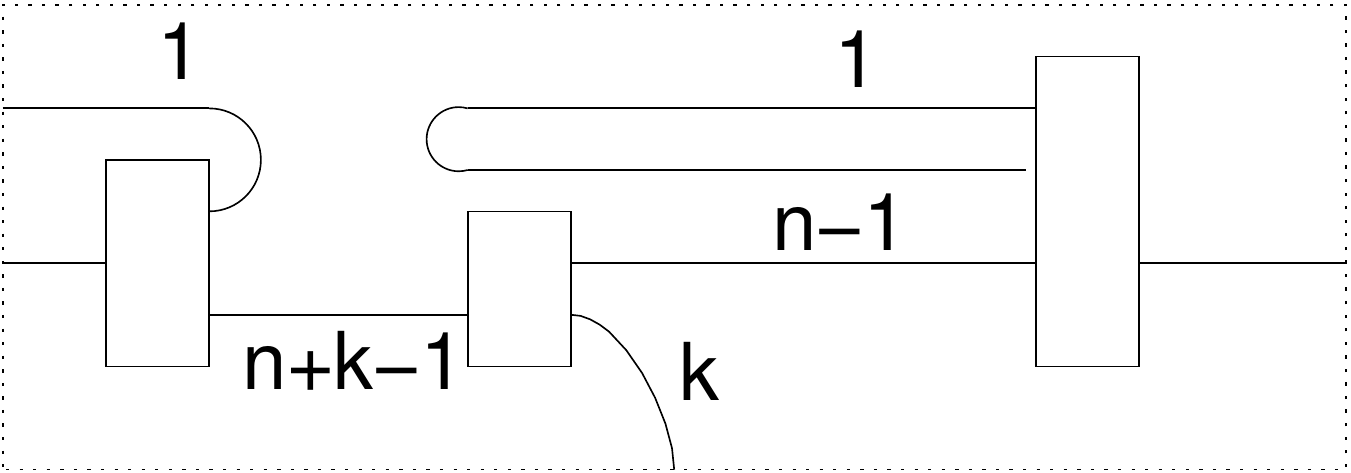} \raisebox{8pt}{$\;- \left(\frac{\Delta_{n+k-2}}{\Delta_{n+k-1}}\right)\;$} \includegraphics[width=1in]{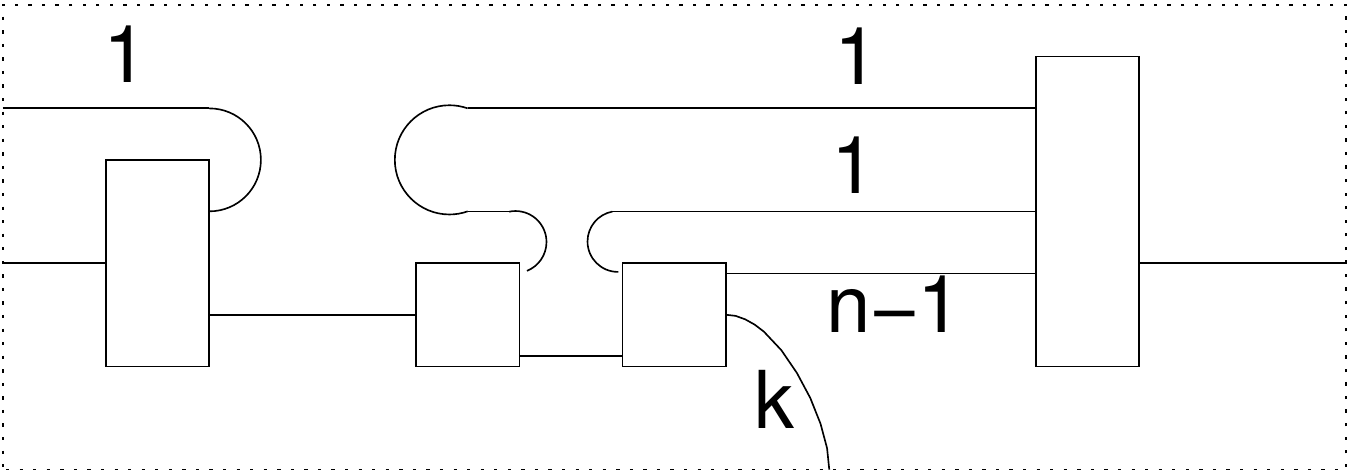}\\
&\raisebox{8pt}{$=$}& \raisebox{8pt}{$- \left(\frac{\Delta_{n+k-2}}{\Delta_{n+k-1}}\right)\;$} \includegraphics[width=1in]{Skeinlemma4v2.pdf}\\
&\raisebox{8pt}{$=$}& \raisebox{8pt}{$- \left(\frac{\Delta_{n+k-2}}{\Delta_{n+k-1}}\right)\;$} \includegraphics[width=1in]{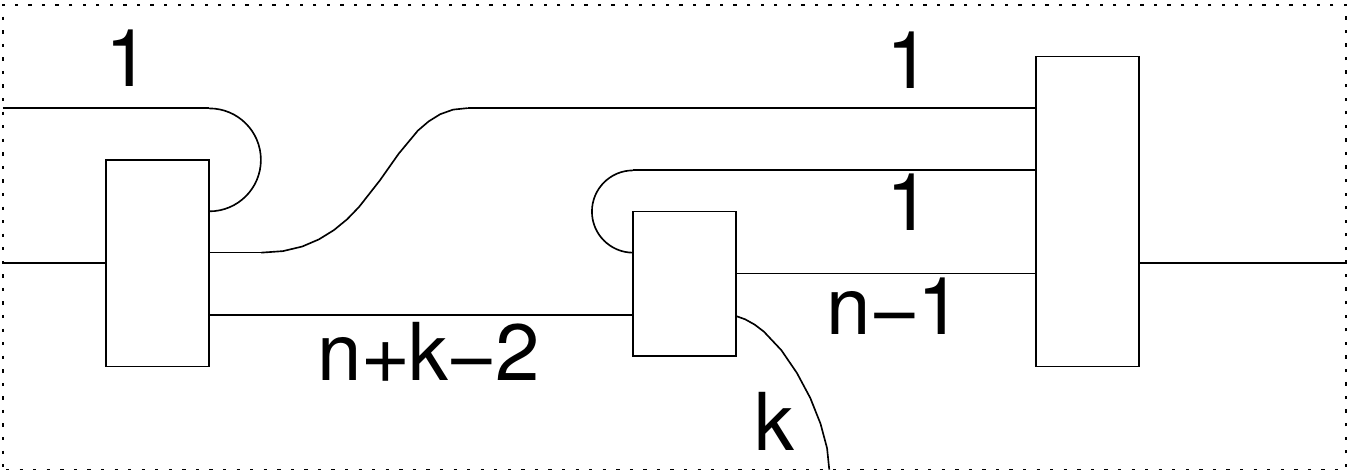}
\end{eqnarray*}

Now when we apply this recursive formula again, the first term will again be zero, and we can continue this process until we get.
$$\includegraphics[width=1in]{Skeinlemma3v2.pdf} \raisebox{9pt}{$\;= - \left(\frac{\Delta_{k-1}}{\Delta_{n+k-1}}\right)\;$} \includegraphics[width=1in]{Skeinlemma6v2.pdf}$$
\end{proof}

Consider a circle $s$ in $S_B$. The circle $s$ appears in $S^{(n+1)}_B$, although it runs through several idempotents. The goal of the argument is to remove one copy of the circle $s$ from the idempotents. Once this is done for each circle in $S_B$, then $S^{(n+1)}_B$ will have been reduced to $S^{(n)}_B$.

Because the diagram $D$ is alternating, the circle $s$ bounds a disk which does not contain any of the other circles in $S_B$. This means that in $S^{(n+1)}_B$, the circle $s$ looks like Figure \ref{loc}. Here all of the arcs ar labeled $n+1$.

\begin{figure}[htbp]
   \centering
   \includegraphics[width=2.0in]{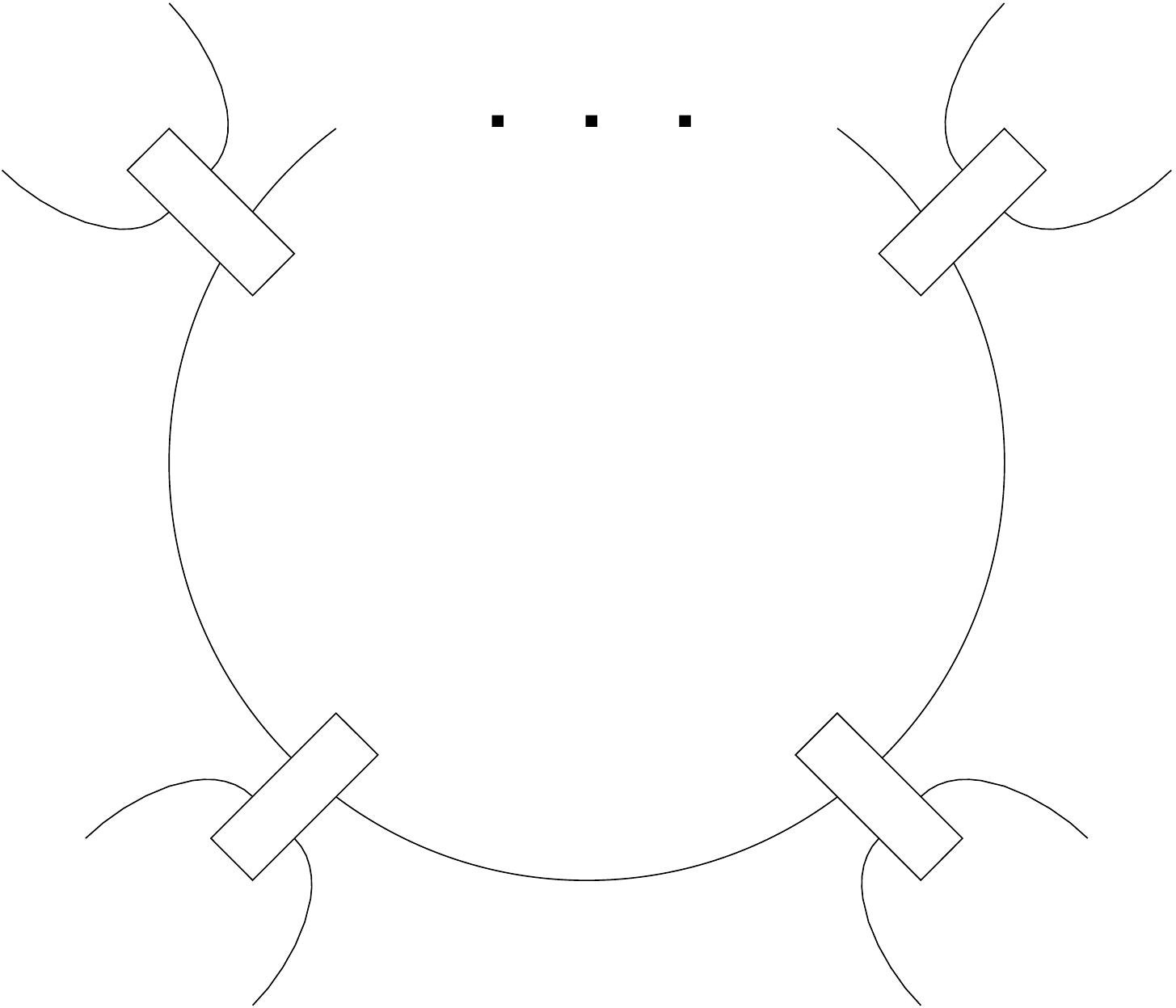} 
   \caption{A circle from $S_B$ seen in $S^{(n+1)}_B$}
   \label{loc}
\end{figure}

 Apply Lemma \ref{skeinreduc} to get the following relation:
$$\includegraphics[width=1.5in]{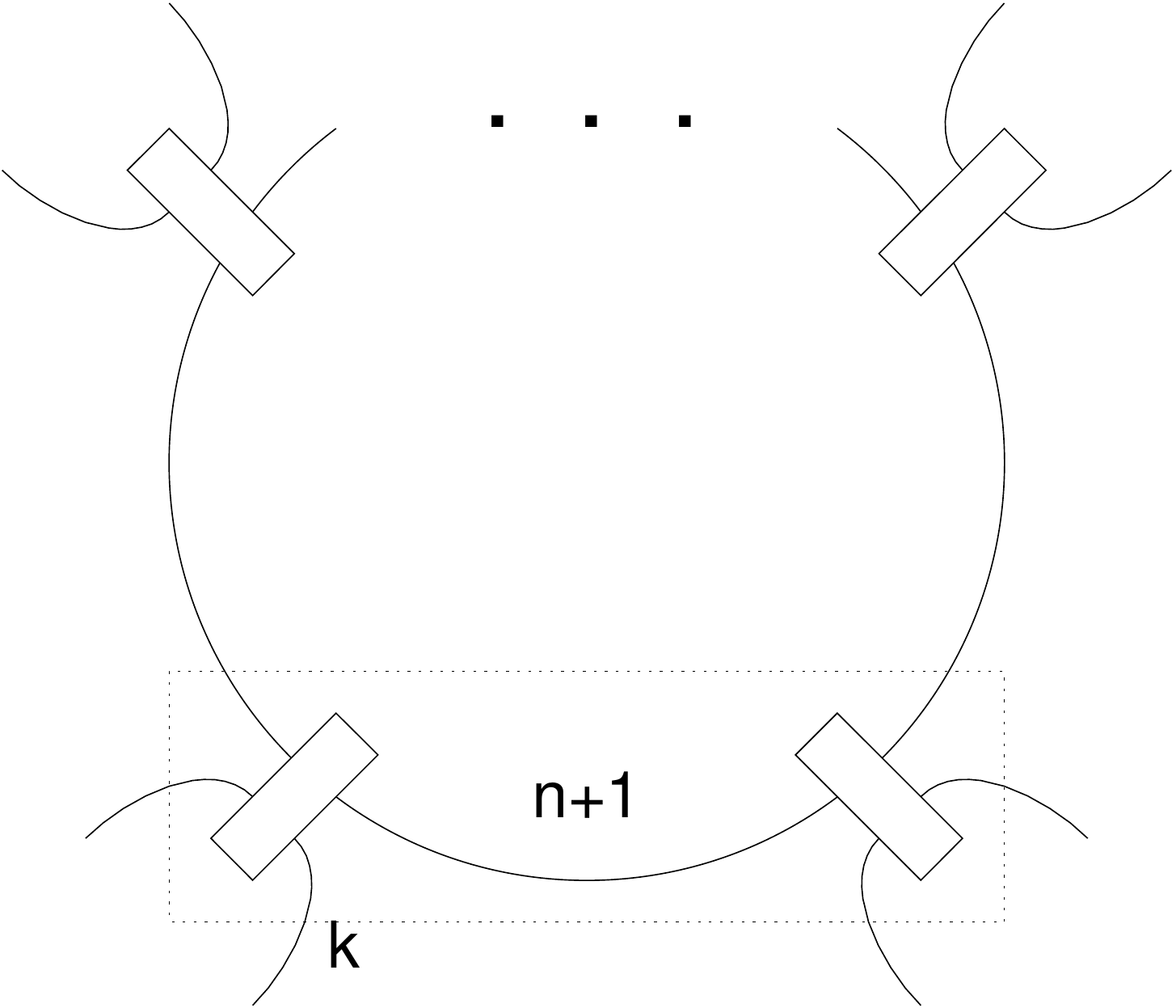} \raisebox{42pt}{$\;=\;$} \includegraphics[width=1.5in]{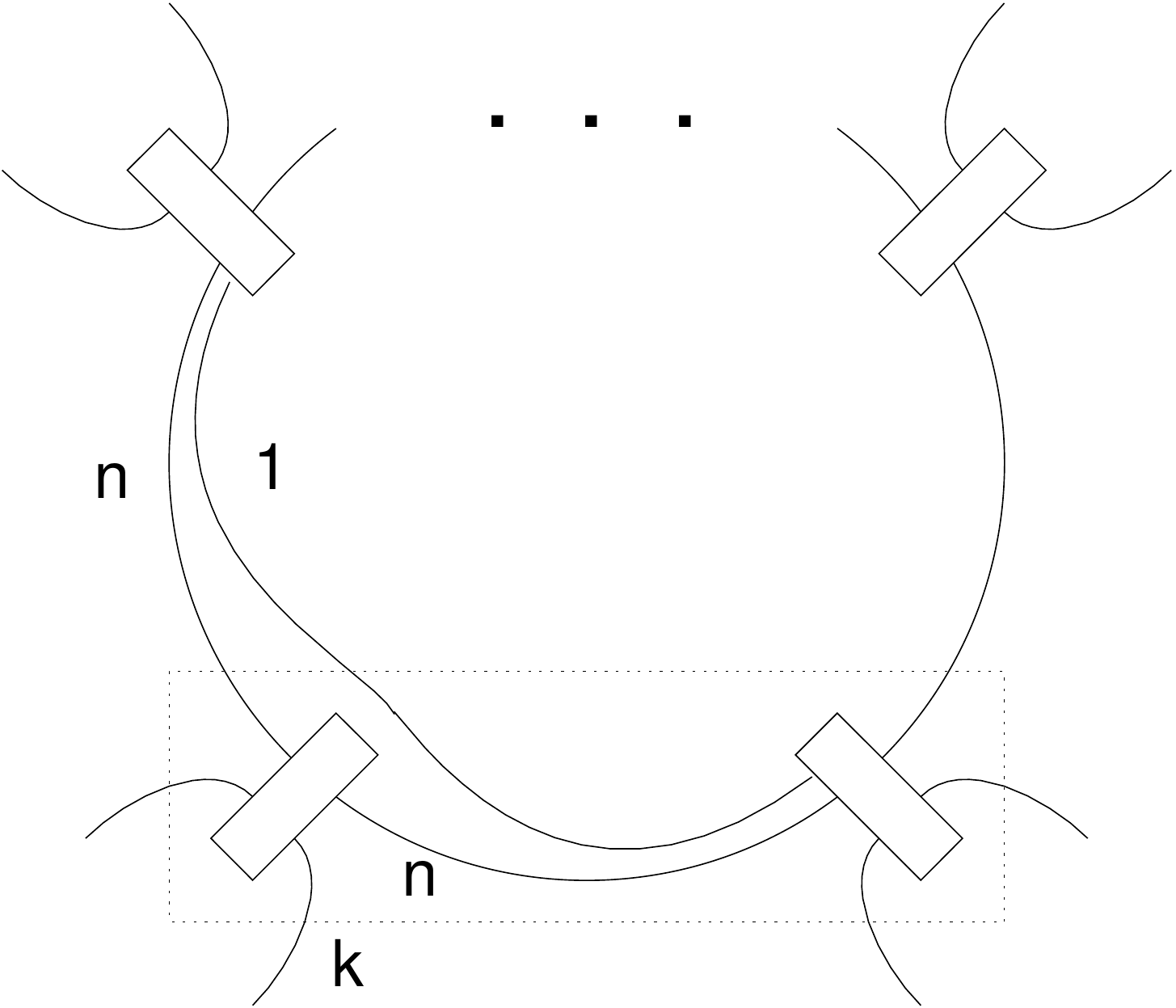} \raisebox{42pt}{$\;- \left(\frac{\Delta_{k-1}}{\Delta_{n+k}}\right)\;$} \includegraphics[width=1.5in]{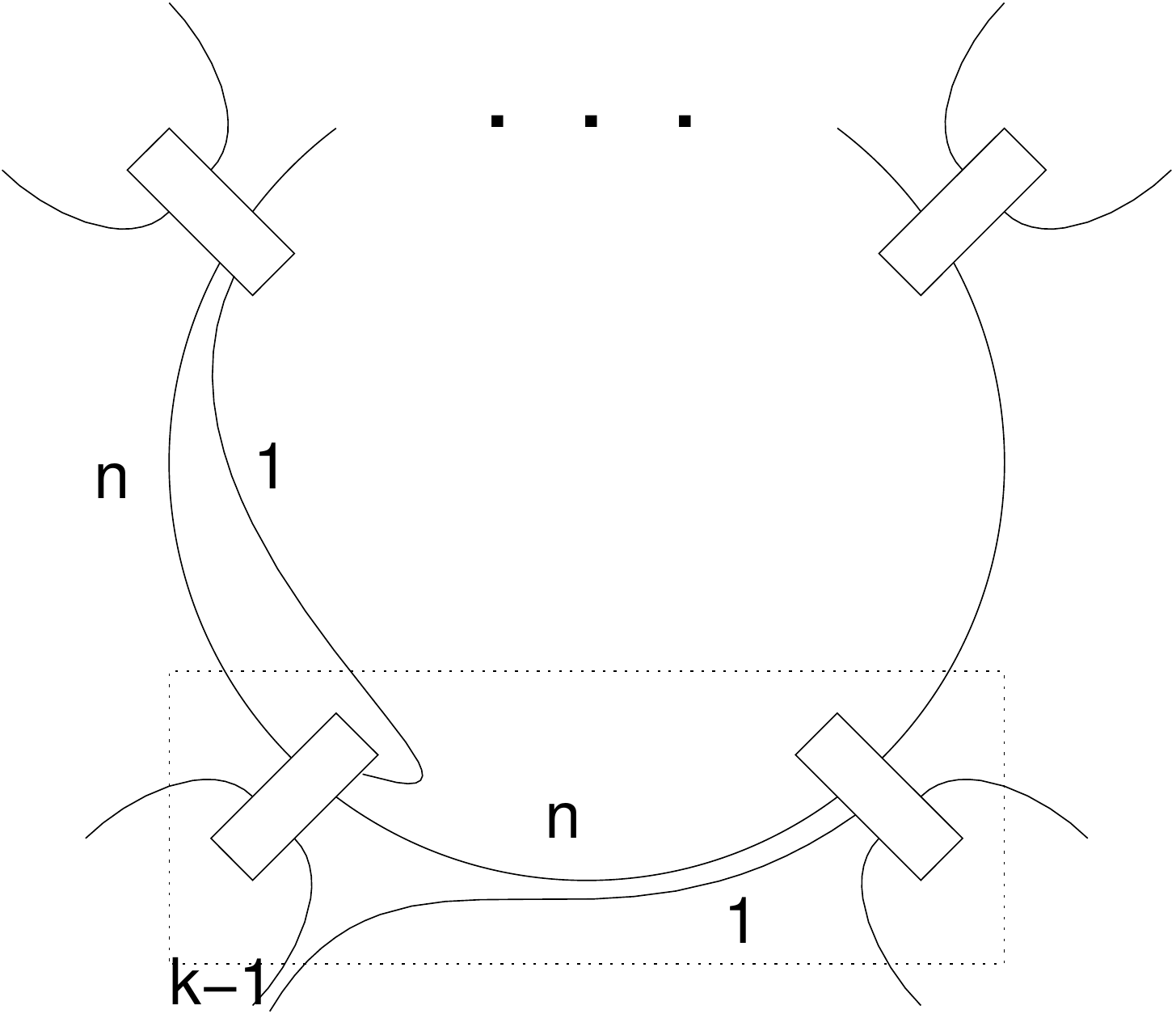}$$

This argument will be applied to each circle in succession, so $k$ is either $n$ or $n+1$ depending on whether the argument has been applied to that circle yet. All non-labelled arcs are either $n$ or $n+1$.

Now $S^{(n+1)}_B$ is expressed as the sum of two terms, and the claim is that the minimum degree of the second term is at least the minimum degree of the first term plus $4(n+1)$. Thus the equation simplifies to:
$$\includegraphics[width=1.5in]{SBlocalSkein1.pdf} \raisebox{42pt}{$\;\dot{=}_{4(n+1)}\;$} \includegraphics[width=1.5in]{SBlocalSkein2.pdf}$$

To see that the claim is true, first note that
$$d\left(\frac{\Delta_{k-1}}{\Delta_{n+k}}\right) = 2(n+1).$$
Now we need to compare the degree of the two diagrams involved. By Lemma \ref{replace} we can get a lower bound for the minimum degrees of these diagrams, and since the knot diagram that these came from was B-adequate, the first term will be an adequate diagram, and thus, the lower bound will be equal to the actual minimum degree. Note that the element in $TL_{n+2}$ shown in Figure \ref{hhh} can be expressed as $h_{n+1}h_n\ldots h_1$. In the previous equation, this element appears in the right most term. When comparing this terms to the first term, each $h_i$ merges two circles into one circle. Thus the number of circles in the diagrams differ by $n+1$. And, finally, a circle can be removed and replaced with a factor of $-A^2-A^{-2}$. This tells us that the difference in the minimum degrees of the diagrams is at least $2(n+1)$. Putting this together with the difference in degrees of the coefficients, the difference in minimum degrees of the terms themselves is at least $4(n+1)$.

\begin{figure}[htbp] %
   \centering
   $\includegraphics[height=1.5in]{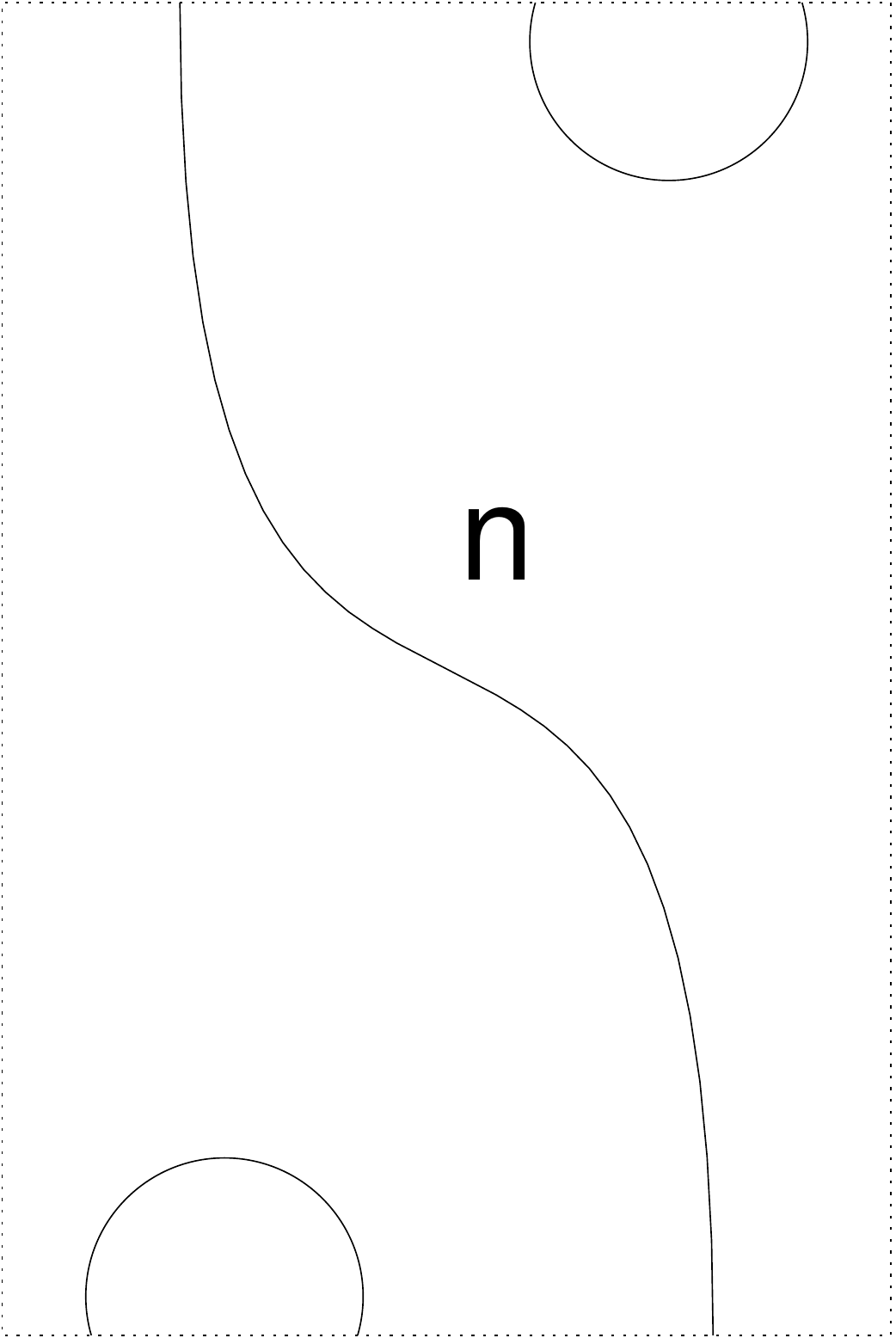} \raisebox{45pt}{$\;=\;$} \includegraphics[height=1.5in]{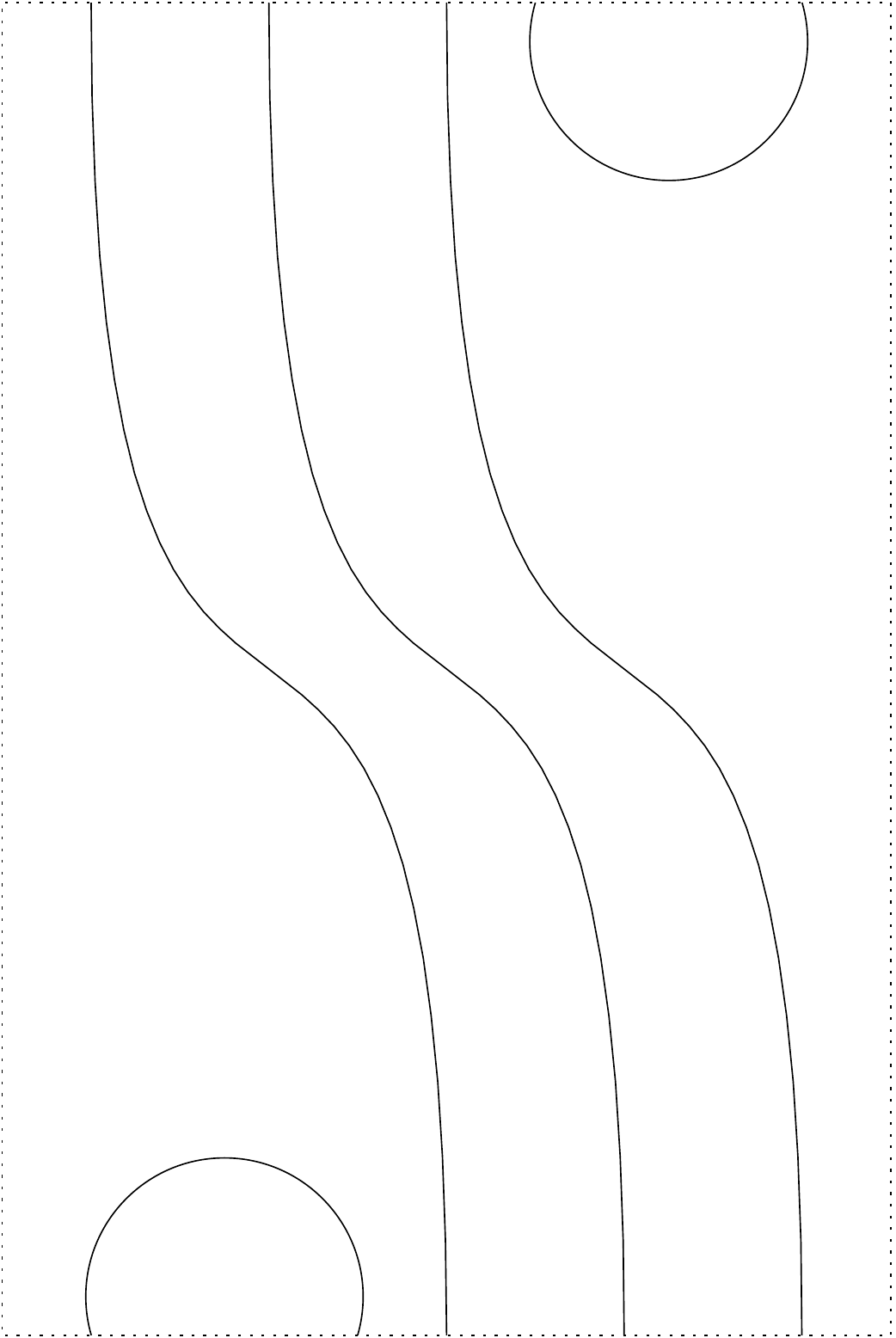} \raisebox{45pt}{$\;= \;$}\includegraphics[height=1.5in]{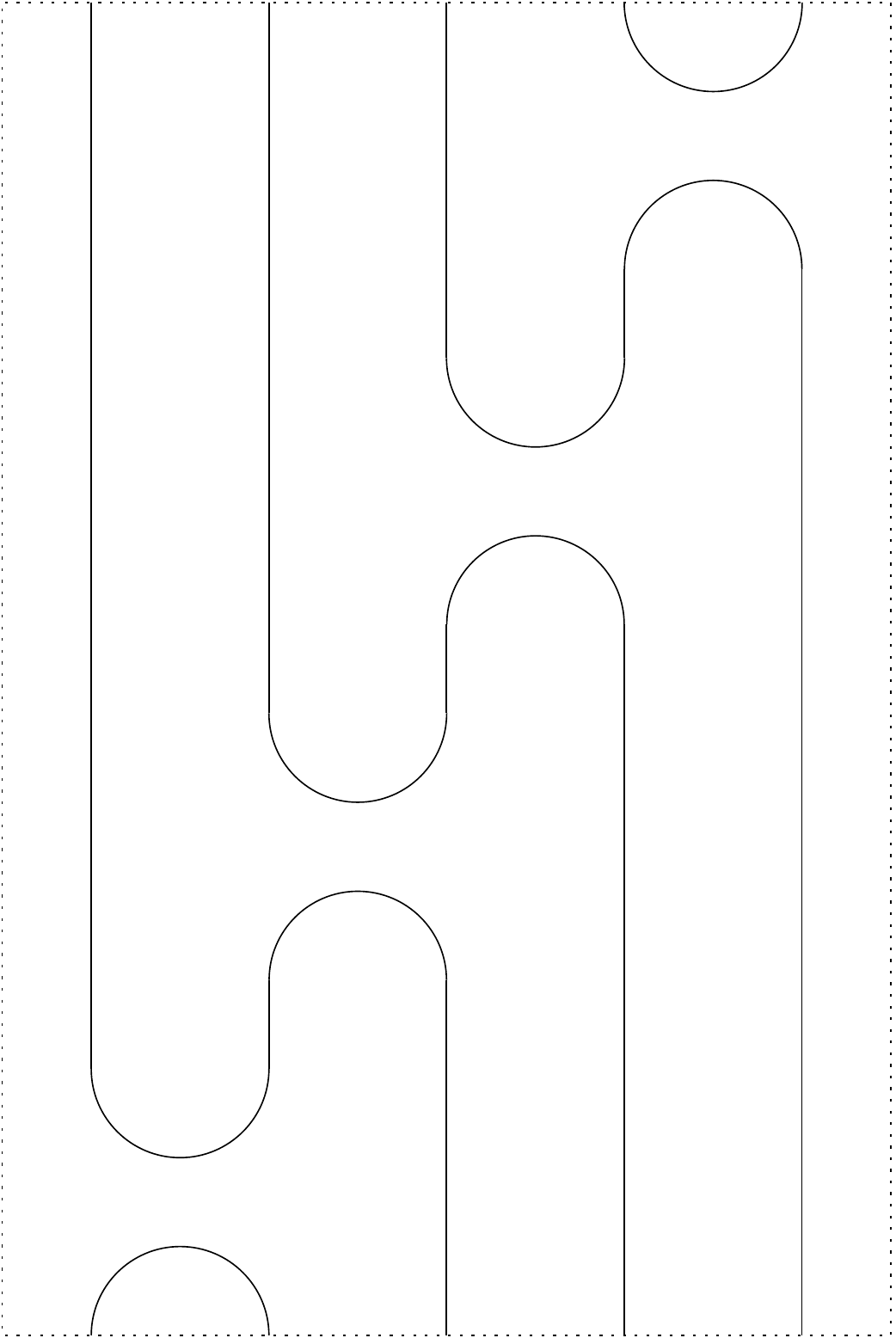}$
   \caption{Multiple pictures expressing $h_{n+1}h_n\ldots h_1$}
   \label{hhh}
\end{figure}

Apply this argument around the circle up to the final idempotent connected to that circle. Now the diagram looks like Figure \ref{final}.

\begin{figure}[htbp] %
   \centering
   $\includegraphics[width=1.5in]{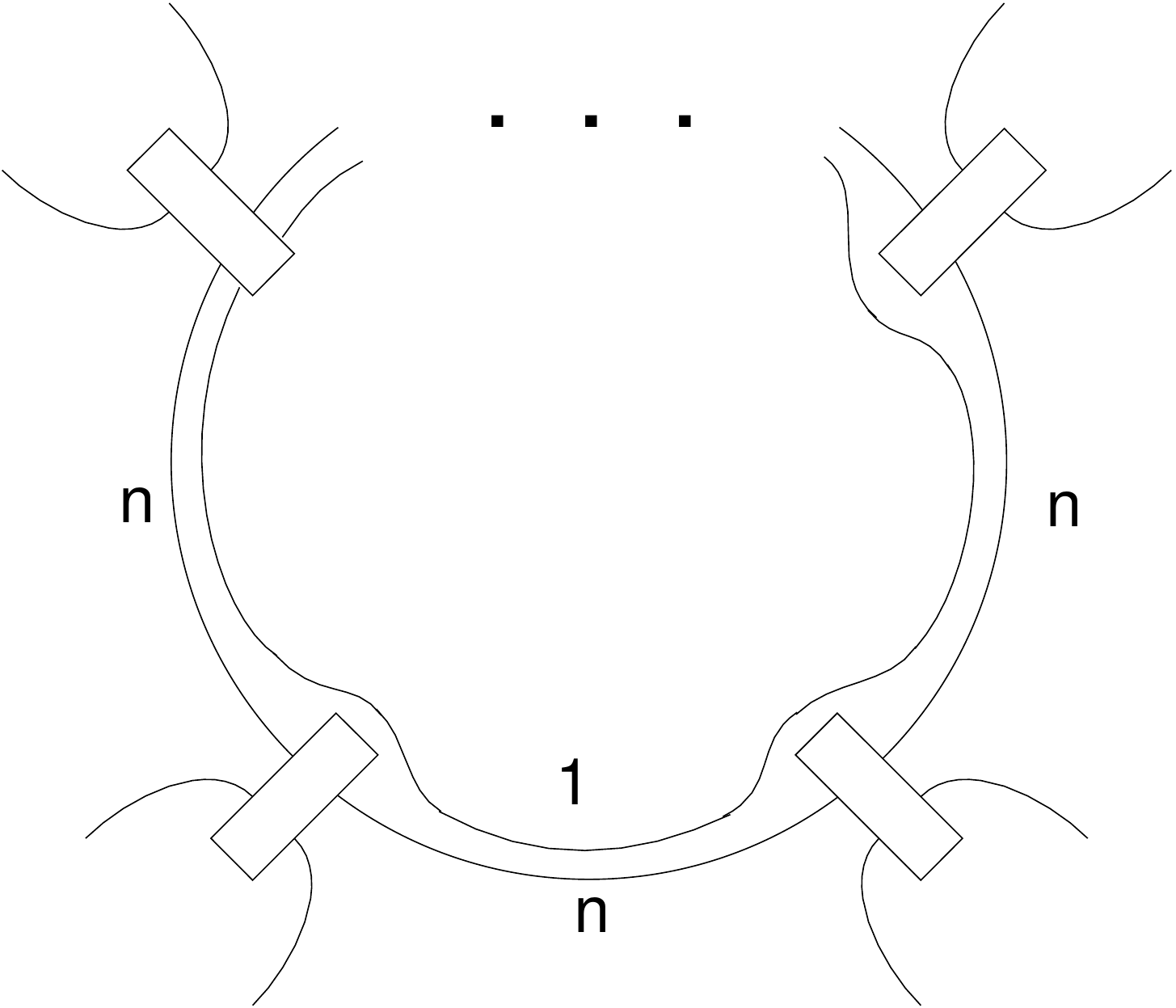} \raisebox{42pt}{$\;=\;$} \includegraphics[width=1.5in]{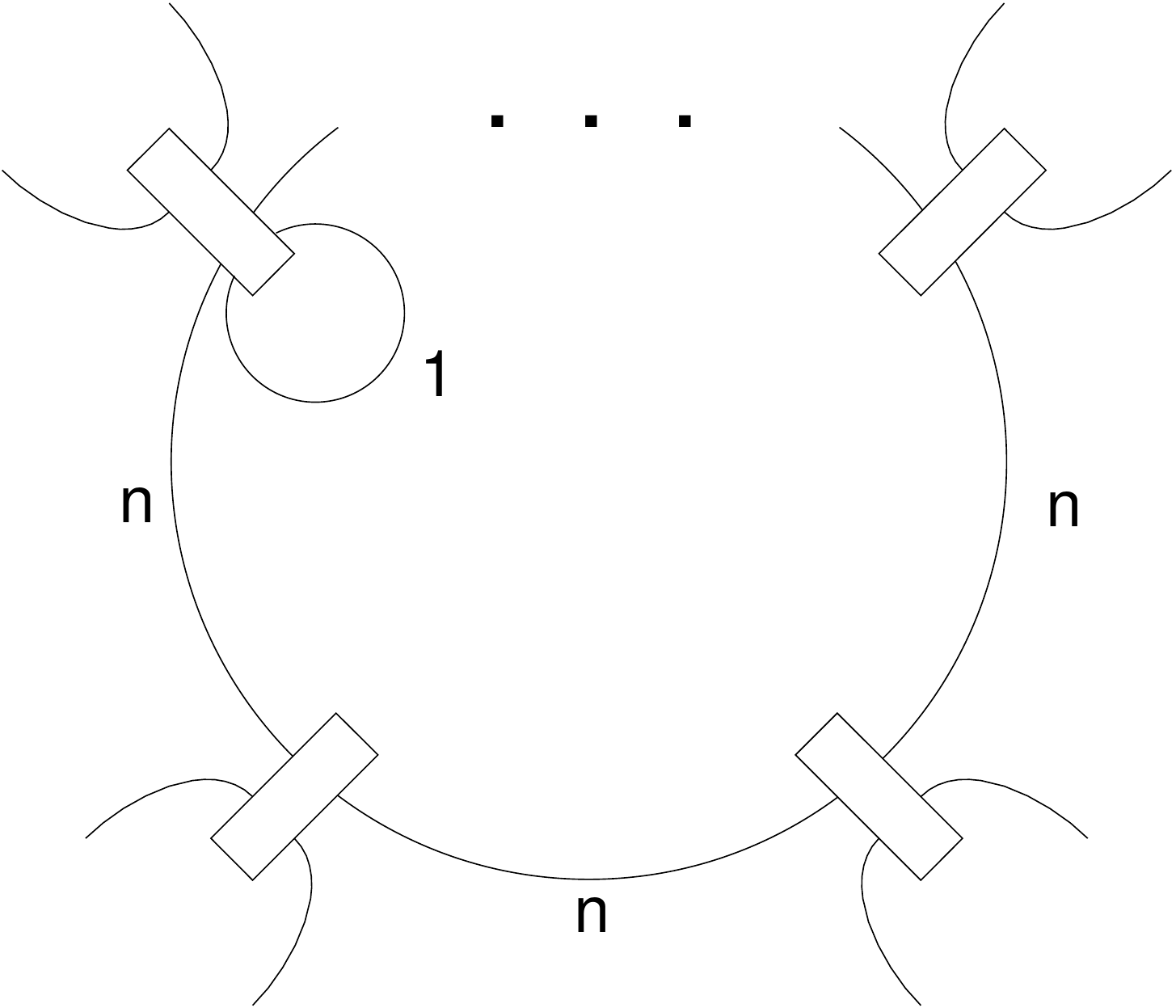} \raisebox{42pt}{$\;= \left(\frac{\Delta_{n+k+1}}{\Delta_{n+k}}\right)\;$}\includegraphics[width=1.5in]{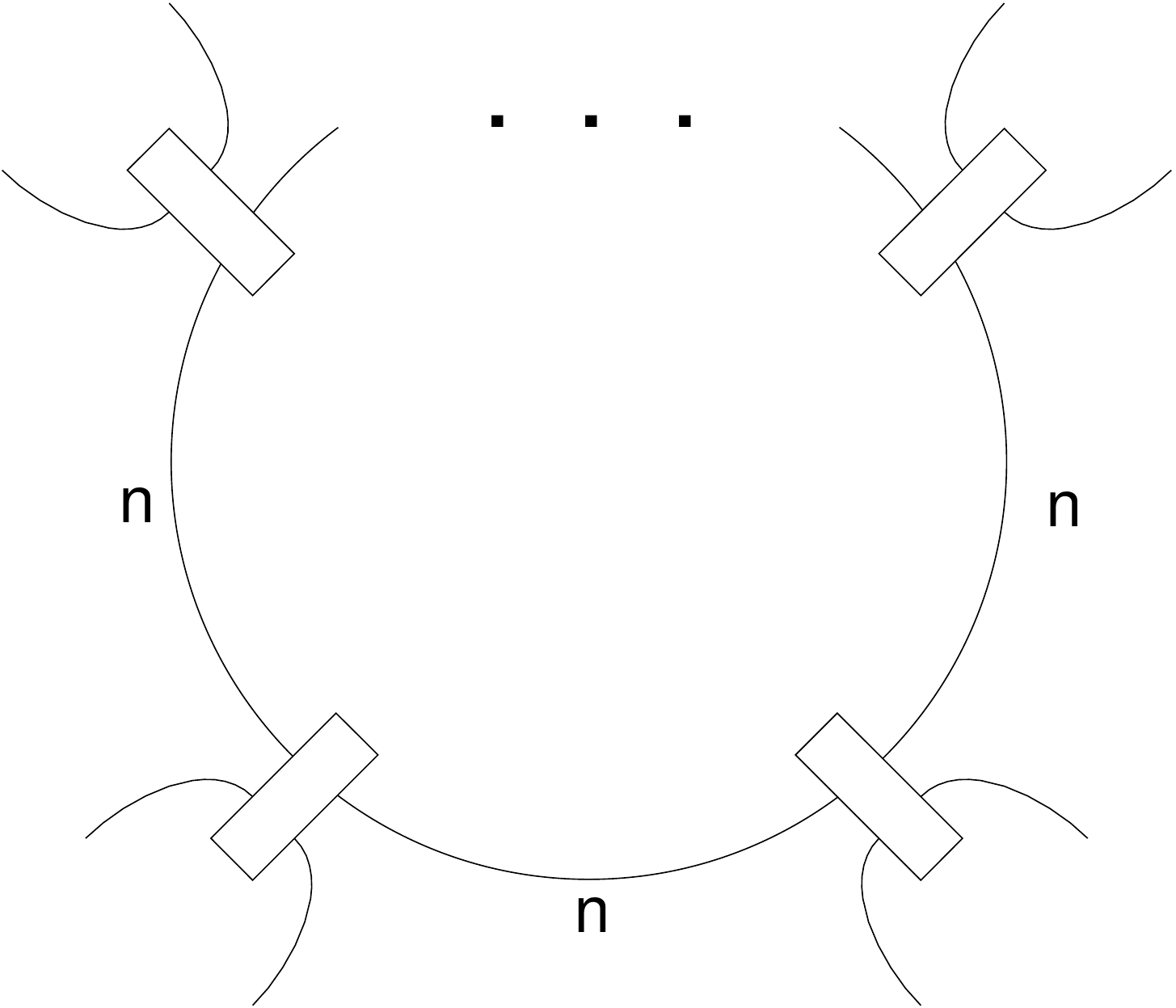}$
   \caption{Reducing $S^{(n+1)}_B$ to $S^{(n)}_B$}
   \label{final}
\end{figure}

We can rewrite the coefficient here:
$$\left(\frac{\Delta_{n+k+1}}{\Delta_{n+k}}\right) = (-1)\frac{A^{2(n+k+2)}-A^{-2(n+k+2)}}{A^{2(n+k+1)}-A^{-2(n+k+1)}} = (-A^{-2})\frac{A^{4(n+k+2)}-1}{{A^{4(n+k+1)}-1}} \;\dot{=}_{4(n+1)}\; 1$$.

Now applying this argument to every circle in $S^{(n+1)}_B$, we see that 
$$S^{(n+1)}_B \;\dot{=}_{4(n+1)}\; S^{(n)}_B$$
and thus,
$$\Gamma_n \;\dot{=}_{4(n+1)}\; \Gamma_{n+1}.$$

This proves Theorem \ref{restated} in the case of reduced alternating diagrams. For the case when the diagram $D$ is B-adequate, most of the proof still applies. The only thing that goes wrong is that Figure \ref{loc} is not accurate because a circle $s$ in $S_B$ might not bound a disk, and thus in $S^{(n+1)}$, may have idempotents which alternate which side of the circle it fuses to other circles. Figure \ref{adtref} shows a non-alternating B-adequate diagram of the trefoil where the dotted cicle is an example such a circle $s$. We would still like to pull out one copy of $s$, but in this case we must be more careful while doing so.

\begin{figure}[htbp] %
   \centering
   \includegraphics[width=1.5in]{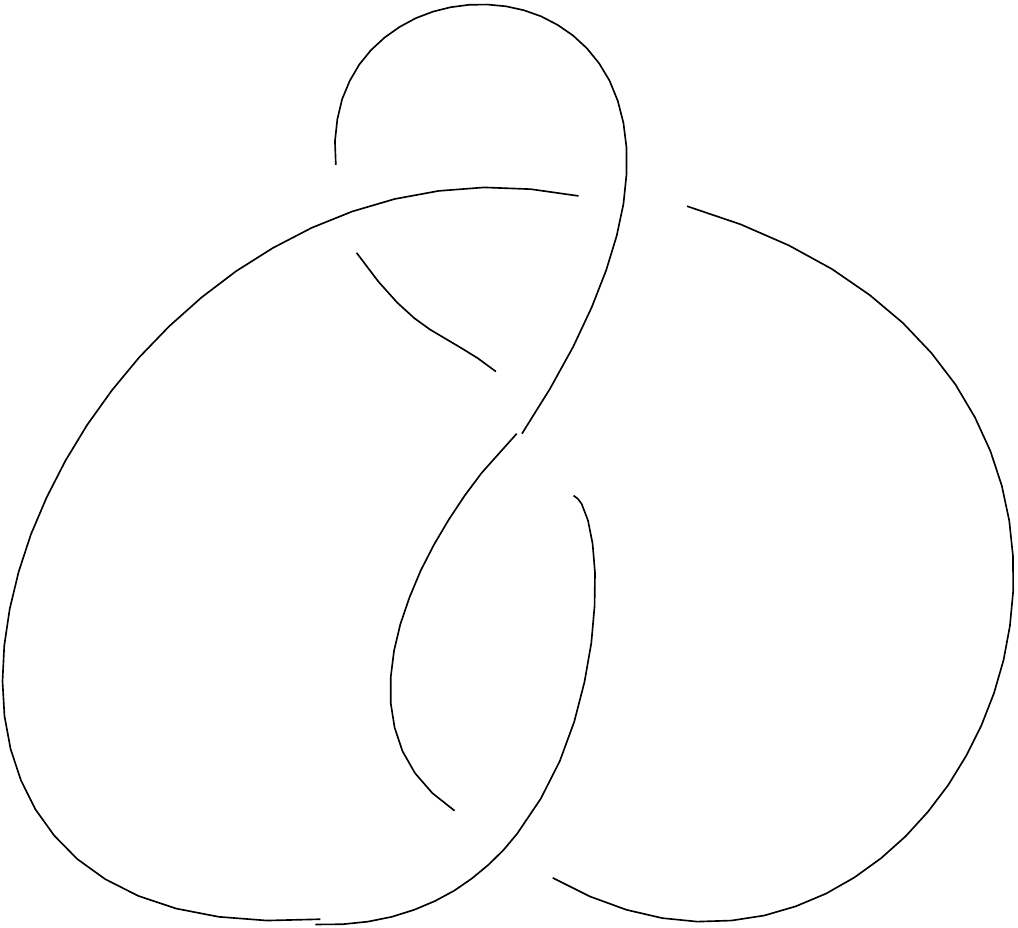}
   \raisebox{42pt}{$\longrightarrow$}
   \includegraphics[width=1.5in]{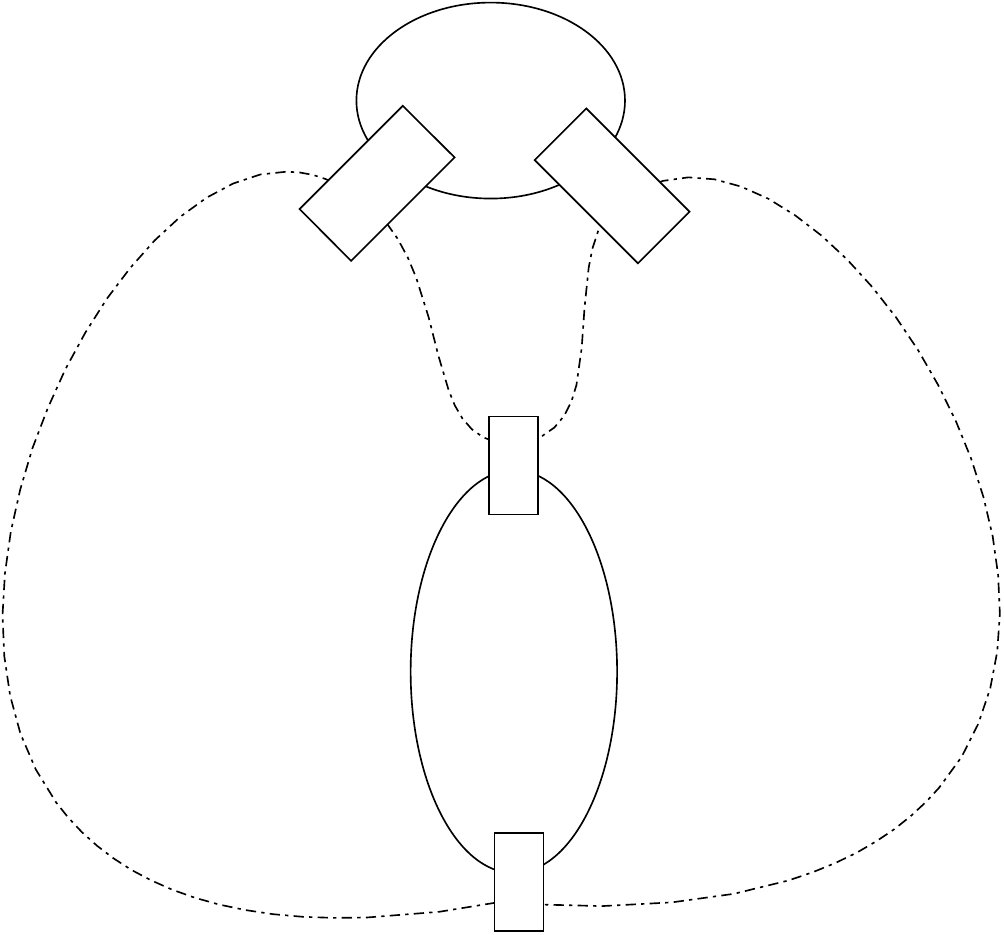} 
   \caption{Example of a non-alternating B-adequate knot diagram}
   \label{adtref}
\end{figure}

First we will modify the diagram as in Figure \ref{twisting} by adding crossings along the circle $s$ between any pair of idempotents which alternate wihch side of $s$ is the outer side. The procedure is to modify the diagram so that the outer strand passes over all of the other copies of $s$, so that it is still the outer strand when it meets the next idempotent in line. Call this new diagram $T$. When expanding $T$ by summing over all possible smoothings of the crossings, only one state is non-zero, and that state is $S^{(n+1)}$. Since this particular smoothing has an equal number of A and B smoothings, we get that $S^{(n+1)} = T$.

\begin{figure}[htbp] %
   \centering
   \includegraphics[height=.8in]{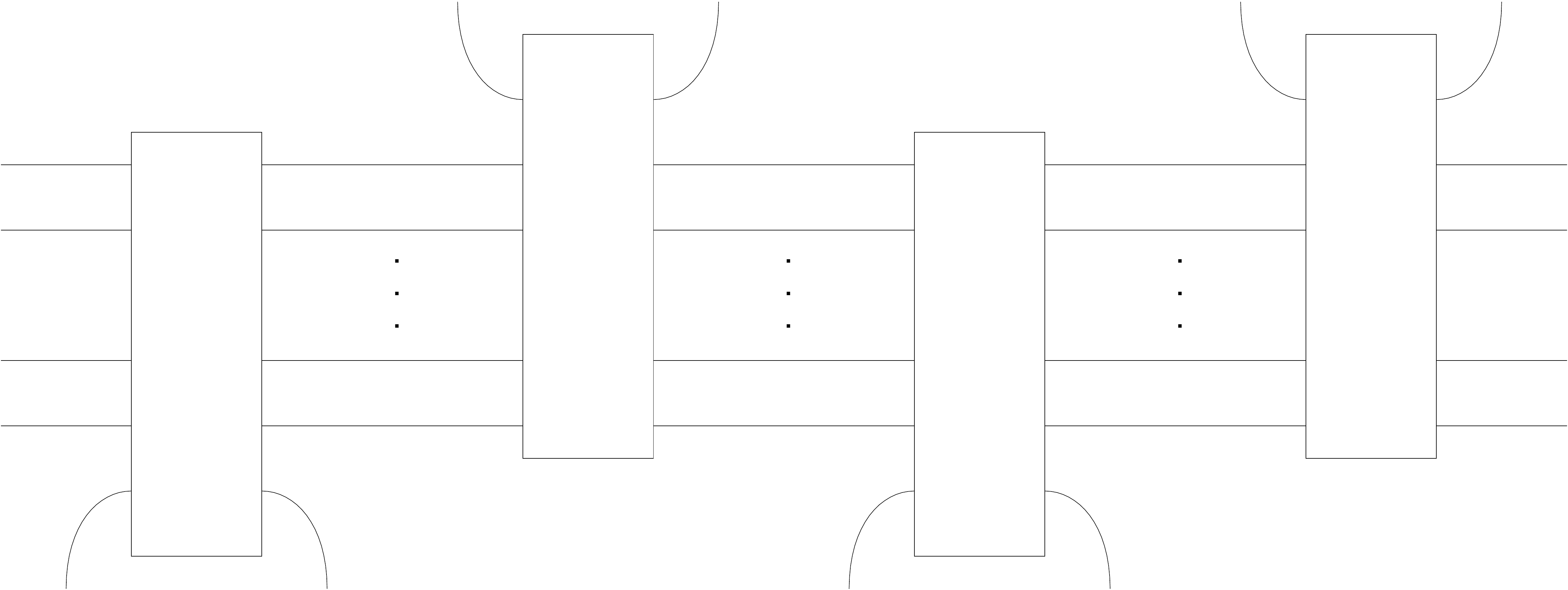}
   \raisebox{24pt}{$\longrightarrow$}
   \includegraphics[height=.8in]{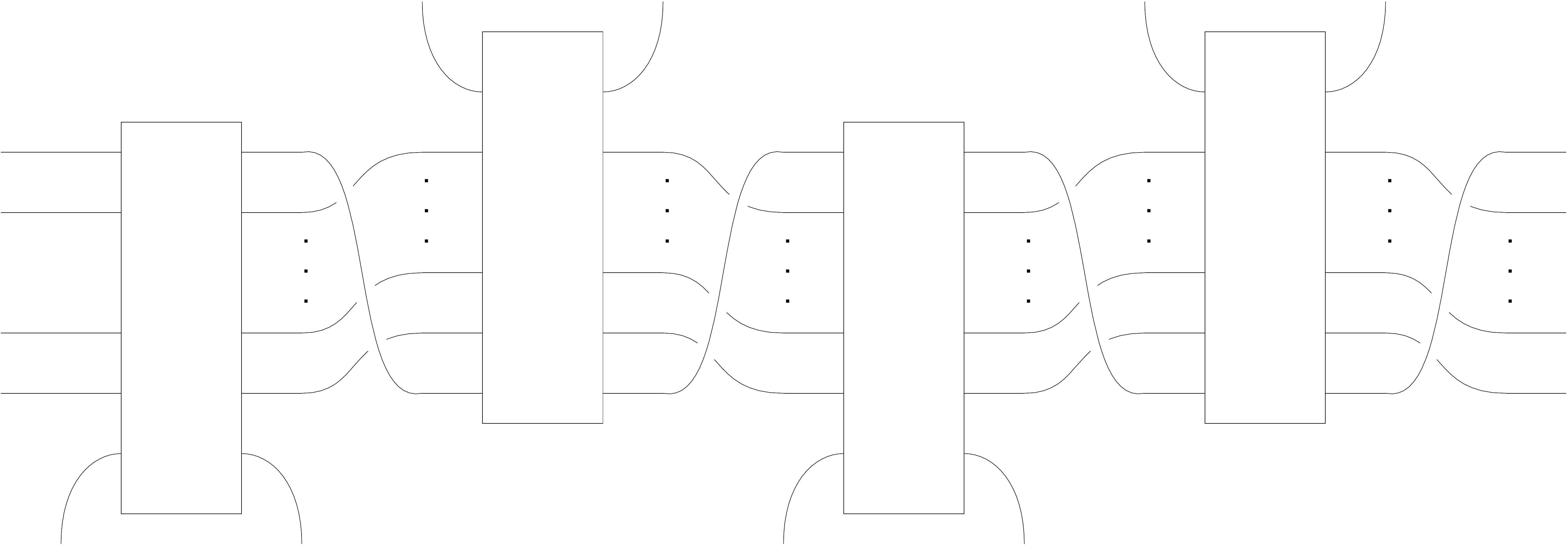}
   \caption{Modifying $S^{(n+1)}$ to get $T$}
   \label{twisting}
\end{figure}

Now we will still apply the procedure to try to pull out one copy of $s$ as in the alternating case. The argument still applies with no modification, as long as we note that in $\bar{T}$, there are still $n+1$ parallel copies of $s$, and they are still unknotted because one of the copies is commpletely over the other copies. Thus that copy can be straightened out to lie next to the other copies. When applying Lemmma \ref{skeinreduc}, the second term is still has $n+1$ fewer circles than $\bar{T}$ because of the same argument. Finally we have the equation:
\begin{eqnarray*}
\includegraphics[height=.8in]{BadS.pdf} & \raisebox{24pt}{$=$}& \includegraphics[height=.8in]{BadT.pdf}\\
&\raisebox{24pt}{$=$}& \includegraphics[height=.8in]{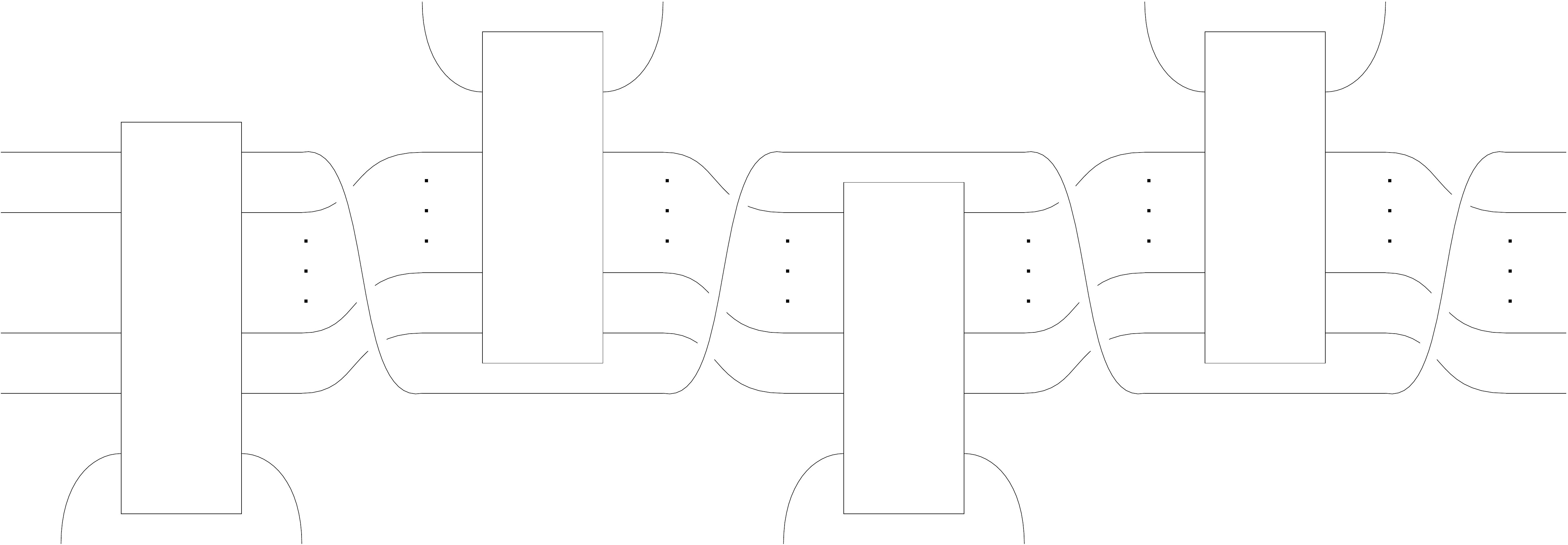}\\
&\raisebox{24pt}{$=$}& \includegraphics[height=.8in]{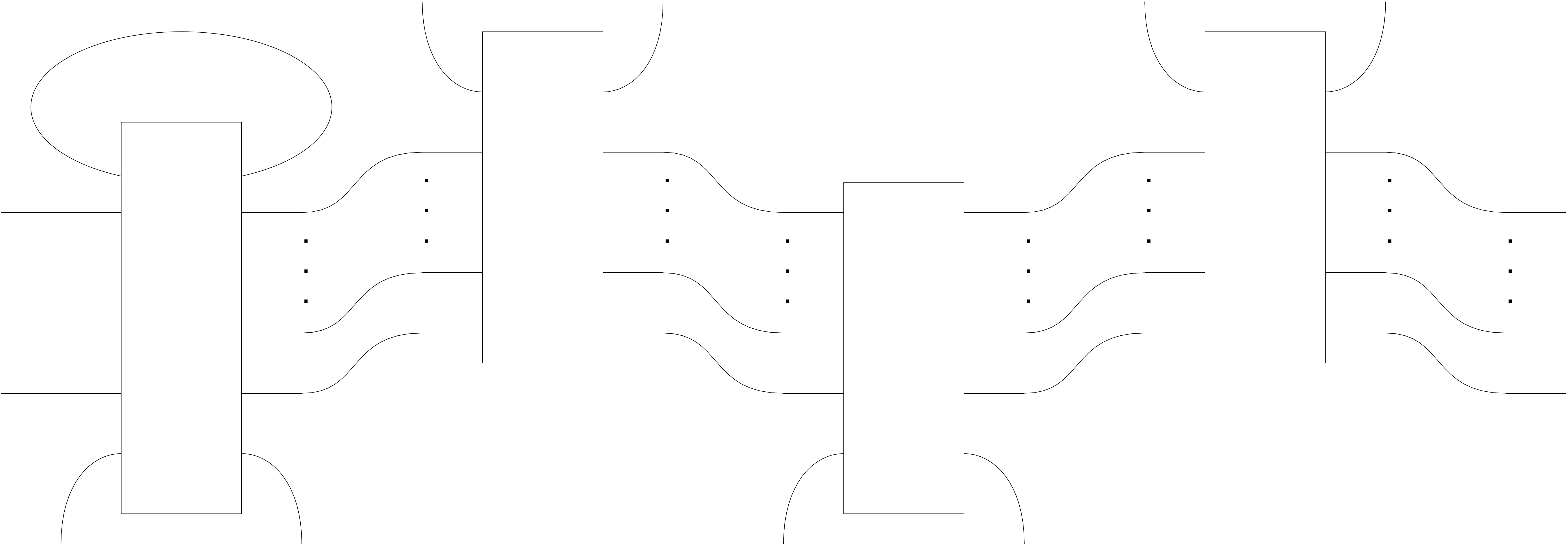}\\
&\raisebox{22pt}{$=$}& \raisebox{22pt}{$\left(\frac{\Delta_{n+k+1}}{\Delta_{n+k}}\right)$} \includegraphics[height=.7in]{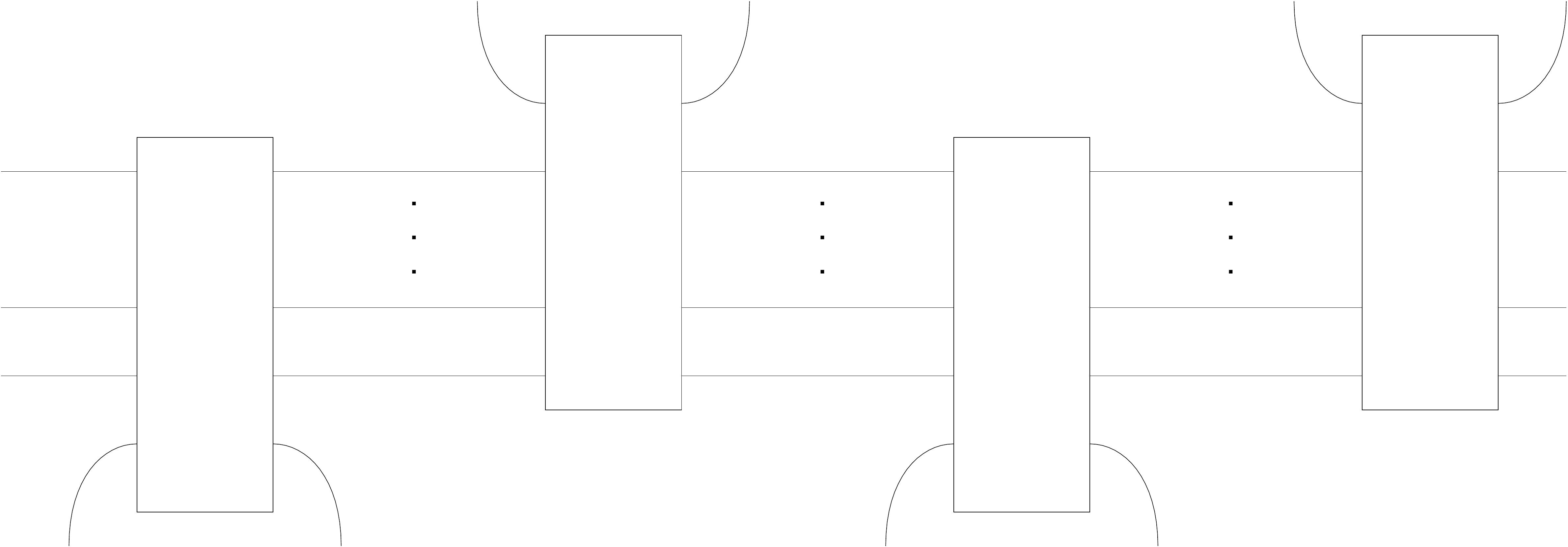}\\
&\raisebox{22pt}{$\dot{=}_{4(n+1)}$}& \includegraphics[height=.7in]{BadTdonev2.pdf}\\
\end{eqnarray*}
which completes the argument.

\end{proof}

\end{document}